\documentclass[]{siamart220329}
\usepackage[utf8]{inputenc}
\usepackage{makecell}
\usepackage{changes}

\headers{Krylov-aware stochastic trace estimation}{T. Chen and E. Hallman}

\title{Krylov-aware stochastic trace estimation\thanks{
\textbf{Funding:} This material is based on work supported by the National Science Foundation under grants DGE-1762114 and DMS-1745654. Any opinions, findings, and conclusions or recommendations expressed in this material are those of the authors and do not necessarily reflect the views of the National Science Foundation.}}
\author{Tyler Chen\thanks{New York University, \href{mailto:tyler.chen@nyu.edu}{\texttt{tyler.chen@nyu.edu}}}
\and 
Eric Hallman\thanks{North Carolina State University, \href{mailto:eric.r.hallman@gmail.com}{\texttt{eric.r.hallman@gmail.com}}}}
\date{}

\usepackage{algorithm}
\usepackage{algorithmicx}
\usepackage[noend]{algpseudocode}
\usepackage{multirow}
\usepackage{amsmath,amssymb}
\usepackage{mathtools}

\usepackage{array}

\usepackage{booktabs}
\usepackage{enumitem}
\usepackage[caption=false]{subfig} 

\usepackage{upgreek}
\usepackage{bm}

\newcommand{\tr}{\operatorname{tr}}
\renewcommand{\vec}{\mathbf}
\newcommand{\T}{{\mkern-1mu\mathsf{T}}}
\newcommand{\F}{\mathsf{F}}
\newtheorem{remark}{Remark}

\usepackage{relsize}
\newcommand{\hpp}{\text{Hutch\raisebox{0.35ex}{\relscale{0.75}++}}}

\usepackage{graphicx}

\newcommand{\EE}{\mathbb{E}}
\newcommand{\VV}{\mathbb{V}}
\newcommand{\PP}{\mathbb{P}}

\usepackage[capitalize]{cleveref}

\usepackage{xcolor}

\begin{document}

\maketitle

\begin{abstract}
    We introduce an algorithm for estimating the trace of a matrix function $f(\vec{A})$ using implicit products with a symmetric matrix $\vec{A}$.
    Existing methods for implicit trace estimation of a matrix function tend to treat matrix-vector products with $f(\vec{A})$ as a black-box to be computed by a Krylov subspace method.
    Like other recent algorithms for implicit trace estimation, our approach is based on a combination of deflation and stochastic trace estimation.
    However, we take a closer look at how products with $f(\vec{A})$ are integrated into these approaches which enables several efficiencies not present in previously studied methods.
    In particular, we describe a Krylov subspace method for computing a low-rank approximation of a matrix function by a computationally efficient projection onto Krylov subspace.
\end{abstract}


\section{Introduction}

We consider the task of estimating  $\tr(f(\vec{A})) = \sum_{i=1}^{d} f(\lambda_i)$ under the assumption that $\vec{A}$ is accessed by matrix-vector products (matvecs). Here $f:\mathbb{R}\to\mathbb{R}$ is a scalar function, $\vec{A}\in \mathbb{R}^{d\times d}$ is a symmetric matrix with eigenvalues $\{\lambda_i\}$, and $f(\vec{A})$ is the corresponding matrix function.
Many algorithms for this problem can be broken into two components: 
(i) computing the trace of an arbitrary symmetric matrix $\vec{B}\in\mathbb{R}^{d\times d}$ accessing $\vec{B}$ only by matvecs; and
(ii) approximating matvecs with $\vec{B} = f(\vec{A})$ by means of a Krylov subspace method.
Each of these individual tasks has been studied extensively, and several analyses have aimed to balance the costs of the two components \cite{persson2022improved,han2017approximating,ubaru2017applications,chen2022randomized}.

It is well known that $\EE[\bm{\uppsi}^\T {\vec{B}} \bm{\uppsi}] = \tr({\vec{B}})$ if $\bm{\uppsi}\in\mathbb{R}^{d}$ satisfies $\EE[\bm{\uppsi}\bm{\uppsi}^\T] = \vec{I}$.
The quadratic trace estimator\footnote{Quadratic trace estimators are sometimes called Hutchinson's trace estimators although they were used prior to Hutchinson's paper \cite{hutchinson1989stochastic}.} $\bm{\uppsi}^\T \vec{B} \bm{\uppsi}$ forms the backbone of the most common stochastic trace estimation algorithms.
When the entries of $\bm{\uppsi}$ are independent and identically distributed (iid) standard Gaussians, it is known that the estimator has variance $2\|\vec{B}\|_\F^2$, where $\|\cdot\|_\F$ denotes the Frobenius norm. 
Thus, if the spectrum of $\vec{B}$ decays quickly, it can be advantageous to compute a low-rank approximation to $\vec{B}$ and apply the quadratic trace estimator to the remainder.
Similar intuition holds for other common choices of $\bm{\uppsi}$ such as iid Rademacher ($\pm 1$) entries or iid Gaussian entries normalized so that $\|\bm{\uppsi}\|_2=\sqrt{d}$.

A number of past works have aimed to combine low-rank approximation with quadratic trace estimation \cite{girard1987algorithme,weisse2006kernel,wu2016estimating,gambhir2017deflation,lin2017randomized,morita2020finite,meyer2021hutch,persson2022improved,baston2022stochastic,epperly2023xtrace}.
Perhaps the most well known is the \hpp{} algorithm \cite{meyer2021hutch} for approximating the trace of an arbitrary implicit matrix\footnote{By ``implicit matrix'', we mean that $\vec{B}$ is assumed to be accessible only through matvecs.} $\vec{B}$.
However, in the case $\vec{B} = f(\vec{A})$, \hpp{} does not take advantage of knowledge that matvecs with $\vec{B} = f(\vec{A})$ are typically approximated by a Krylov subspace method.
Moreover, at least in its original form, \hpp{} and related variants must be run separately for each implicit matrix.
This is in contrast to simple quadratic trace estimation based Krylov subspace methods which essentially produce a quadrature approximation that can be used to simultaneously approximate $\tr(f(\vec{A}))$ for multiple functions $f$ efficiently \cite{bai1996largescale,bai1996bounds,schnack2020accuracy,chen2022randomized}.

The primary goal of this paper is to show that when matvecs with $f(\vec{A})$ are computed via a Krylov subspace method, more efficient algorithms for estimating $\tr(f(\vec{A}))$ are possible.
This is enabled by two critical observations. 
First, rather than treating $f(\vec{A})$ as an arbitrary matrix, it makes sense to think about how the spectrum of $f(\vec{A})$ depends on the spectrum of $\vec{A}$. 
Indeed, the natural primitive operation in our setting is matvecs with $\vec{A}$.
Second, matvecs with $f(\vec{A})$ and a set of vectors can be approximated very efficiently when the vectors themselves are elements of a Krylov subspace generated with $\vec{A}$.

\paragraph{Contributions} Our main contribution is \cref{alg:main}, which combines low-rank approximation with quadratic trace estimation to estimate $\tr(f(\vec{A}))$. It is similar in style to \hpp{} but exploits the structure of a block Krylov space to reduce the required number of matvecs.
In particular, we show how a low-rank approximation to $f(\vec{A})$ can be computed more efficiently, including for multiple functions $f$ simultaneously.
In \cref{sec:greedy} we present the details of this algorithm and discuss its relation to past work. 
In \cref{sec:variants} we present two variants: \cref{alg:adaptive}, an adaptive version that takes as input an error tolerance and failure probability, and \cref{alg:restart}, a version designed for situations with a limited amount of memory available. Numerical experiments in \cref{sec:experiments} show that our methods compare favorably to existing ones. 

\subsection{Motivating example}

Consider the matrix function $f(\vec{A}) = \exp(-\beta \vec{A})$ parameterized by the scalar $\beta>0$.
This function arises in equilibrium quantum thermodynamics as the partition function $Z(\beta) \equiv \tr(\exp(- \beta \vec{A}))$ and gives us access to properties of a quantum system such as the specific heat, magnetization, and entropy.
Evaluating the dependence of $Z(\beta)$ on the inverse Boltzmann temperature $\beta$ is of general interest in the study of quantum systems \cite{weisse2006kernel,schnack2020accuracy,faridfar2021thermodynamic,chen2022randomized}.

\begin{figure}
    \centering
    \includegraphics[width=\textwidth]{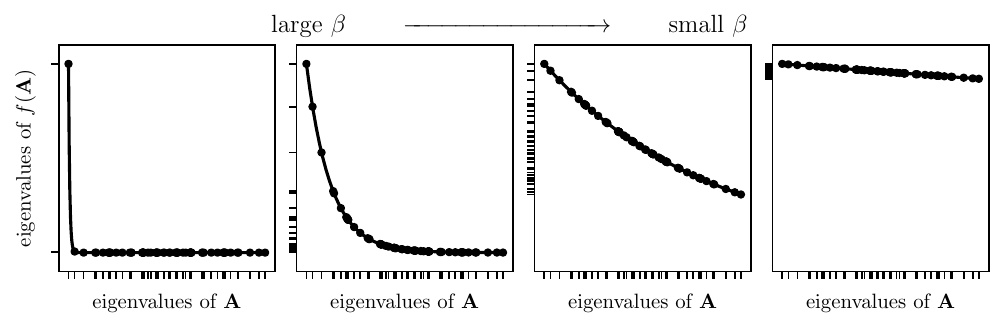}
    \caption{Dependence of eigenvalues of $f(\vec{A}) = \exp(-\beta \vec{A})$ on the parameter $\beta$. 
    Observe that the decay in the spectrum of $\exp(-\beta \vec{A})$ can be extremely fast if $\beta$ is large (left) or extremely slow if $\beta$ is small (right).
    We would like an algorithm which can approximate $\tr(\exp(-\beta \vec{A}))$ for many values of $\beta$ simultaneously.
    }
    \label{fig:exp_example}
\end{figure}

As depicted in \cref{fig:exp_example}, the spectrum of $\exp(-\beta \vec{A})$ depends strongly on the value of $\beta$.
The spectrum decay influences which algorithms are most suited for computing the trace of the matrix function.
For instance, when $\beta$ is large, the partition function $Z(\beta)$ is determined almost entirely by the smallest eigenvalue of $\vec{A}$ (ground state energy) and approaches based on low rank approximation are effective \cite{saibaba2017randomized,morita2020finite,li2021randomized}.
On the other hand, when $\beta$ is small, the spectrum of $\exp(-\beta\vec{A})$ is relatively flat and standard quadratic trace estimators \cite{girard1987algorithme,skilling1989eigenvalues,hutchinson1989stochastic} work very well.
For intermediate $\beta$, a combination of these approaches is effective: the contribution to the trace of the top several eigenvalues can be computed directly using a low-rank approximation, and the contribution of the remaining eigenvalues can be computed using quadratic trace estimators.

\subsection{Notation}
Bolded capital letters $\vec{A}$, $\vec{\Omega}$ denote matrices. Bolded lowercase Roman and Greek letters $\vec{x}$, $\vec{y}$, $\boldsymbol{\uppsi}$ denote vectors. The vector $\vec{e}_i$ denotes the $i$th column of the identity matrix $\vec{I}_n$, whose size $n$ can be inferred from context. Similarly, the matrix $\vec{E}_i$ denotes the Kronecker product $\vec{e}_i\otimes \vec{I}_b$, i.e., columns $(i-1)b+1$ through $ib$ of $\vec{I}_{nb}$, where $n$ and $b$ can be inferred from context. The all-zero matrix is $\vec{0}$, and its dimensions can be inferred from context. 

The transpose of a matrix is $\vec{A}^\T$. 
The number of columns of $\vec{A}$ is $\text{cols}(\vec{A})$.
An \emph{orthonormal} matrix $\vec{Q}$ is one with orthonormal columns: i.e., $\vec{Q}^\T\vec{Q} = \vec{I}$. 
If a matrix has the eigenvalue decomposition $\vec{A} = \vec{U}\vec{\Lambda}\vec{U}^\T$, the matrix function $f(\vec{A})$ is $\vec{U}f(\vec{\Lambda})\vec{U}^\T$, where $f(\vec{\Lambda})$ applies $f$ to each diagonal entry of $\vec{\Lambda}$. For a symmetric matrix $\vec{A}$, a \emph{dominant eigenspace} of $\vec{A}$ of dimension $k$ is one corresponding to $k$ eigenvalues of $\vec{A}$ of largest magnitude. If $\vec{A}$ has no repeated eigenvalues, then the dominant eigenspace is unique for each $k$.
We will write $\lambda_{\textup{max}}$ and $\lambda_{\textup{min}}$ for the largest and smallest eigenvalues of $\vec{A}$, $\|\vec{A}\|_2$ for the operator norm, and $\|\vec{A}\|_\F$ for the Frobenius norm.
We denote by $\sigma_{\textup{max}}(\vec{B})$ and $\sigma_{\textup{min}}(\vec{B})$ the largest and smallest singular values of $\vec{B}$.

The submatrix consisting of rows \( r \) through \( r' \) and columns \( c \) through \( c' \) is denoted by \( [\vec{B}]_{r:r',c:c'} \). A colon with no indices refers to all rows or columns.
Thus, \( [\vec{B}]_{:,1:2} \) denotes the first two columns of \( \vec{B} \), and \( [\vec{B}]_{3,:} \) denotes the third row of \( \vec{B} \).

\section{Background}

\subsection{Implicit trace estimation} \label{bgrd:trace}

Given an orthonormal matrix $\vec{Q}\in \mathbb{R}^{d\times b}$, define
\begin{align*}
    \widehat{\vec{B}} &\equiv \vec{Q}\vec{Q}^\T\vec{B} + (\vec{I}-\vec{Q}\vec{Q}^\T)\vec{B}\vec{Q}\vec{Q}^\T, \\
    \widetilde{\vec{B}} &\equiv (\vec{I}-\vec{Q}\vec{Q}^\T)\vec{B}(\vec{I}-\vec{Q}\vec{Q}^\T),
\end{align*}
and note that the linearity of the trace operator implies that 
\begin{equation}    \label{eqn:trace_decomp}
    \tr(\vec{B}) = \tr(\widehat{\vec{B}}) + \tr(\widetilde{\vec{B}}).
\end{equation}
The cyclic property of the trace implies that 
\[
\tr(\widehat{\vec{B}})
= \tr(\vec{Q}^\T \vec{B} \vec{Q}),
\]
which allows us to efficiently compute this term exactly.
If we define 
$\vec{y}_i \equiv (\vec{I}-\vec{Q}\vec{Q}^\T)\bm{\uppsi}_i$
where $\bm{\uppsi}_i$ is drawn from some spherically symmetric distribution, then the $\vec{y}_i / \| \vec{y}_i\|_2$ are distributed uniformly on the unit hypersphere in the complement of the column span of $\vec{Q}$.
Therefore, 
\[
    (d-b) \EE\left[ \frac{\vec{y}_i\vec{y}_i^\T }{\vec{y}_i^\T \vec{y}_i} \right]
    =  (\vec{I} - \vec{Q} \vec{Q}^\T),
\]
where $(\vec{I} - \vec{Q} \vec{Q}^\T)$ acts as the identity operator on the complement of the column span of $\vec{Q}$.
Again using the cylic property of the trace and that $(\vec{I} - \vec{Q} \vec{Q}^\T)^2 = (\vec{I} - \vec{Q} \vec{Q}^\T)$, we see that
\[
    \tr(\widetilde{\vec{B}})
    =  \tr((\vec{I} - \vec{Q} \vec{Q}^\T) \vec{B} (\vec{I} - \vec{Q} \vec{Q}^\T))
    = (d-b) \EE\left[ \frac{\vec{y}_i^\T \vec{B} \vec{y}_i}{\vec{y}_i^\T \vec{y}_i} \right]
\]
Thus, we obtain an unbiased estimator
\begin{equation}
    \label{eqn:defl_est}
 \tr(\vec{B})
 \approx \tr(\vec{Q}^\T\vec{B}\vec{Q}) + \frac{d-b}{m} \sum_{i=1}^{m} \frac{\vec{y}_i^\T \vec{B} \vec{y}_i}{\vec{y}_i^\T \vec{y}_i}.
\end{equation}

A simple choice of spherically symmetric distribution for $\bm{\uppsi}_i$ is for the entries to be iid Gaussians. 
In this case, while $\vec{y}_i^\T \vec{y}_i$ concentrates around $d-b$ when $d-b$ is large, for finite $d-b$, using the normalizing factors $(d-b) / \vec{y}_i^\T \vec{y}_i$ can significantly improve performance over using the un-normalized vectors $\vec{y}_i$ alone.
The effect is particularly noticeable when the spectrum of $\widetilde{\vec{B}}$ is flat.
The use of such a normalization in the context of stochastic trace estimation was suggested in \cite[\S 2.3]{epperly2023xtrace}. 

The critical observation about an estimator like \cref{eqn:defl_est} is that if $\vec{Q}$ is close to the dominant eigenspace of $\vec{B}$, then $\widetilde{\vec{B}}$ may have Frobenius norm significantly smaller than that of $\vec{B}$.
The estimator $ {\vec{y}_i^\T \vec{B} \vec{y}_i}/{\vec{y}_i^\T \vec{y}_i}$ used in \cref{eqn:defl_est} has variance closely related to $\|\widetilde{\vec{B}}\|_{\F}^2$, so the required number of random vectors used may be reduced significantly compared to if quadratic trace estimation was used with $\vec{B}$ directly.

A natural approach to finding a $\vec{Q}$ near the dominant eigenspace of $\vec{B}$ is by \emph{sketching} \cite{halko2011finding}. 
This can be done by taking $\vec{Q}$ to be an orthonormal basis for $\vec{B} \vec{\Omega}$, where $\vec{\Omega}\in\mathbb{R}^{d\times b}$ is a random test matrix.
Thus, we arrive at \cref{alg:generalized_simple} which returns an estimate to $\tr(\vec{B})$ using $2b + m$ matvecs with $\vec{B}$. 
This is essentially the \hpp{} algorithm of \cite{meyer2021hutch}, except that the latter specifies that $b = m$ and that $\vec{\Omega}$ and $\vec{\Psi}$ have iid Rademacher entries, and the normalization scheme of \cite{epperly2023xtrace} was not used.

\begin{algorithm}[ht]
\caption{Implicit stochastic trace estimation \cite{meyer2021hutch}}\label{alg:generalized_simple}
\fontsize{10}{14}\selectfont
\begin{algorithmic}[1]
\Procedure{implicit-trace}{$\vec{B}, b, m$}
\State sample Gaussian matrices $\vec{\Omega} \in\mathbb{R}^{d\times b}$, $\vec{\Psi}\in \mathbb{R}^{d\times m}$
\State $\vec{Q} = \textsc{orth}(\vec{B} \vec{\Omega})$ 
\State  $t_{\text{defl}} = \tr(\vec{Q}^\T\vec{B}\vec{Q})$ \label{line:tdefl}
\State $\vec{Y} = (\vec{I}-\vec{Q}\vec{Q}^\T)\vec{\Psi}$
\Comment{$\vec{Y} = [\vec{y}_1,\ldots,\vec{y}_{m}]$}
\State $t_{\text{rem}} = \frac{d-b}{m} \sum_{i=1}^{m} {\vec{y}_i^\T \vec{B} \vec{y}_i}/{\vec{y}_i^\T \vec{y}_i}$
\State \Return $t_{\text{defl}} + t_{\text{rem}}$
\EndProcedure
\end{algorithmic}
\end{algorithm}

\subsection{Krylov subspace methods for matrix function approximation}
\label{bgrd:krylov}

It is common to use Krylov subspace methods to approximate products with $f(\vec{A})$ \cite{higham2008functions,frommer2008matrix}. 
Specifically, expressions of the form $f(\vec{A})\vec{Z}$ and $\vec{Z}^\T f(\vec{A}) \vec{Z}$ can be approximated using information from the block Krylov subspace
\begin{equation}
\label{eqn:krylov_subspace}
\mathcal{K}_{q+1}(\vec{A},\vec{Z}) 
= \operatorname{span}\{ \vec{Z}, \vec{A}\vec{Z}, \ldots, \vec{A}^q \vec{Z} \}.
 \end{equation}
Here the span is interpreted as the span of the union of the columns of the constituent matrices. 

The block-Lanczos algorithm \cite{golub1977block} takes a starting matrix $\vec{Z}\in \mathbb{R}^{d\times b}$ with QR factorization $\vec{Z}=\vec{Q}_1\vec{R}_1$, and after $q$ steps, computes a matrix
\begin{equation}\label{def:qbar}
    \vec{\bar{Q}}_{q+1} = \left[\vec{Q}_1,\vec{Q}_2,\ldots,\vec{Q}_{q+1}\right] \in \mathbb{R}^{d\times (q+1)b}
\end{equation}
whose columns form an orthonormal basis for $\mathcal{K}_{q+1}(\vec{A},\vec{Z})$.
The block columns of $\vec{\bar{Q}}_{q+1}$ satisfy a symmetric block-tridiagonal recurrence
\begin{equation}
\label{eqn:krylov_recurrence}
    \vec{A} \vec{\bar{Q}}_{q} = \vec{\bar{Q}}_{q} \vec{T}_{q} + \vec{Q}_{q+1} \vec{R}_{q+1} \vec{E}_{q}^\T
\end{equation}
where
\begin{align*}
    \vec{T}_{q} = \begin{bmatrix}
    \vec{M}_1 & \vec{R}_2^\T \\
    \vec{R}_2 & \ddots & \ddots \\
    & \ddots & \ddots & \vec{R}_{q}^\T \\
    &&\vec{R}_{q} & \vec{M}_{q}
    \end{bmatrix} \in \mathbb{R}^{qb\times qb}.
\end{align*}
We call a recurrence of the form \eqref{eqn:krylov_recurrence} a block-Lanczos recurrence.\footnote{The dimensions of $\vec{\bar{Q}}_{q+1}$ and $\vec{T}_q$ are given under the assumption that the blocks $\vec{R}_1,\vec{R}_2,\ldots$ all have rank $b$.}

\Cref{alg:block_lanczos} presents pseudocode for the block-Lanczos algorithm. 
For reasons which will become apparent
, we also include an additional input parameter $n$, so that the algorithm runs for $q+n$ iterations but only performs the reorthogonalization step in \hyperref[line:reorth]{line \ref{line:reorth}} for the first $q$ iterations. 
In order to keep the presentation as simple as possible, we assume for \cref{alg:block_lanczos} that the Krylov subspace $\mathcal{K}_{q+n}(\vec{A},\vec{Z})$ has dimension $(q+n)b$.
If this assumption is not met, the algorithm can terminate unsuccessfully.
A fully practical implementation of the block-Lanczos algorithm should handle rank-deficient blocks by deflation or other means; see for instance \cite{zhou2008block}.

\begin{algorithm}[ht]
\caption{block-Lanczos}\label{alg:block_lanczos}
\fontsize{10}{14}\selectfont
\begin{algorithmic}[1]
\Procedure{block-Lanczos}{$\vec{A}, \vec{Z}, q, n$}
\State \(\vec{Q}_1, \vec{R}_1 = \textsc{qr}(\vec{Z}) \)
\For {\( k=1,2,\ldots,q+n\)}
    \If{$k=1$}
    \State $\vec{Z} = \vec{A}\vec{Q}_1$
    \Else
    \State \( \vec{Z} = \vec{A} \vec{Q}_{k} - \vec{Q}_{k-1} \vec{R}_{k}^\T \)
    \EndIf
    \State \( \vec{M}_k = \vec{Q}_{k}^\T \vec{Z}  \)\label{line:mmat}
    \State \( \vec{Z} = \vec{Z} - \vec{Q}_{k} \vec{M}_k \) \label{line:zmat}
    \If{$2\leq k \leq q$}
    \State \( \vec{Z} = \vec{Z} - \vec{\bar{Q}}_{k-1} (\vec{\bar{Q}}_{k-1}^\T  \vec{Z} ) \) \Comment{reorthogonalize against $\vec{\bar{Q}}_{k-1}$ \eqref{def:qbar}}   \label{line:reorth}
    \EndIf
    \State \(\vec{Q}_{k+1}, \vec{R}_{k+1} = \textsc{qr}(\vec{Z}) \) \label{line:qnext}
\EndFor
\State \Return $\vec{T}_{q+n}$, $\vec{\bar{Q}}_{q+1}$, optionally $\vec{R}_1$
\EndProcedure
\end{algorithmic}
\end{algorithm}

After $q$ steps of the process have been carried out, we obtain the approximation 
\begin{align}
    f(\vec{A})\vec{Z}&\approx \vec{\bar{Q}}_{q} f(\vec{T}_{q})\vec{\bar{Q}}_{q}^\T\vec{Z}
    = \vec{\bar{Q}}_{q}\left[f(\vec{T}_{q})\right]_{:,1:b}\vec{R}_1\label{eqn:mv_approx}\\
    \intertext{and the quadratic form approximation}
    \vec{Z}^\T f(\vec{A})\vec{Z}&\approx \vec{Z}^\T\vec{\bar{Q}}_{q}f(\vec{T}_{q})\vec{\bar{Q}}_{q}^\T\vec{Z}
    = \vec{R}_1^\T \left[f(\vec{T}_{q})\right]_{1:b,1:b}\vec{R}_1,\label{eqn:quad_approx}
\end{align}
where $\vec{Z} = \vec{Q}_1\vec{R}_1$. 
These approximations are exact when $f$ is a low-degree polynomial.
\begin{lemma}\label{thm:poly_exact}
Suppose the Krylov subspace $\mathcal{K}_q(\vec{A},\vec{Z})$ has dimension $qb$ (so that no rank-deficient blocks are encountered). Then the approximation \cref{eqn:mv_approx} is exact when $f$ is a polynomial of degree at most $q-1$, and the approximation \cref{eqn:quad_approx} is exact when $f$ is a polynomial of degree at most $2q-1$.
\end{lemma}
Proofs of all lemmas and theorems are given in \cref{section:proofs}.

\Cref{thm:poly_exact} yields a simple bound on the rate of convergence of these approximations.
\begin{lemma} \label{thm:func_rate}
Under the assumptions of \cref{thm:poly_exact}, for $q>0$, the approximation \cref{eqn:mv_approx} satisfies
\[
    \| f(\vec{A})\vec{Z} - \vec{\bar{Q}}_{q}\left[f(\vec{T}_{q})\right]_{:,{1:b}}\vec{R}_1 \|_2
    \leq 2 \| \vec{Z} \|_2  \min_{\deg(p) < q} \left( \max_{x\in [\lambda_{\textup{min}}, \lambda_{\textup{max}}]} | f(x) - p(x)| \right),
\]
and the approximation \cref{eqn:quad_approx} satisfies
\[
    \| \vec{Z}^\T f(\vec{A})\vec{Z} - \vec{R}_1^\T \left[f(\vec{T}_{q})\right]_{{1:b,1:b}}\vec{R}_1 \|_2
    \leq 2 \| \vec{Z} \|_2^2  \min_{\deg(p) < 2q} \left( \max_{x\in [\lambda_{\textup{min}}, \lambda_{\textup{max}}]} | f(x) - p(x)| \right).
\]
\end{lemma}

Assuming the interval $[\lambda_{\textup{min}},\lambda_{\textup{max}}]$ is known, \cref{thm:func_rate} allows us to apply standard bounds from approximation theory to choose $n$ \cite{trefethen2019approximation}.
For instance, if $\|\vec{A}\|_2\leq 1 $ and $f$ is analytically continuable to a Bernstein ellipse $E_\rho$ where it satisfies  $|f(x)|\leq M$ for all $x\in E_\rho$, then the error of the best polynomial approximation to $f$ on $[-1,1]$ is bounded by $M \rho^{-n} / (\rho-1)$.
Bounds for non-analytic functions, including $\nu$-times differentiable functions whose ($\nu+1$)-th derivative is of bounded variation, can also be found in \cite{trefethen2019approximation}.

The Lanczos algorithm without reorthogonalization is highly susceptible to the impacts of finite precision arithmetic.
While the rates of convergence of the approximations \cref{eqn:mv_approx,eqn:quad_approx} are typically reduced when the algorithm is run in finite precision arithmetic, the bound in \cref{thm:func_rate} can still be expected to hold to close degree, \emph{even without reorthogonalization}.
This has been shown rigorously for the block size one case \cite{druskin1991error,knizhnerman1996simple,musco2018stability}. 
In practice, the bounds in \cref{thm:func_rate} are often very pessimistic, even in finite precision arithmetic. 
Developing spectrum-dependent a priori and a posteriori bounds suitable for use as practical stopping criteria is an active area of research \cite{frommer2008stopping,frommer2009error,ilic2009restarted,frommer20132norm,frommer2014convergence,frommer2015error,chen2022error,xu2022posteriori}.

Alternate methods for matrix function approximation include explicit polynomial approaches such as those based on Chebyshev series \cite{weisse2006kernel,han2017approximating}.
While these approaches tend to satisfy similar bounds to \cref{thm:func_rate}, they typically do not outperform Lanczos-based methods \cite{chen2022randomized,chen2022error}.

\subsection{A simple algorithm for the trace of a matrix function}

As mentioned in the introduction, it is common to use the ideas from \cref{bgrd:krylov} to implement the matvecs with $f(\vec{A})$ in algorithms like \cref{alg:generalized_simple}.
Such an approach has been used in numerical experiments in \cite{meyer2021hutch,persson2022improved}, and an implementation is described explicitly in \cref{alg:generalized_simple_func}.

\begin{algorithm}[ht]
\caption{Implicit stochastic trace estimation for functions}\label{alg:generalized_simple_func}
\fontsize{10}{14}\selectfont
\begin{algorithmic}[1]
\Procedure{implicit-trace-func}{$\vec{A}, f, b, q, m, n$}
\State sample Gaussian matrices $\vec{\Omega} \in\mathbb{R}^{d\times b}$, $\vec{\Psi}\in \mathbb{R}^{d\times m}$
\State $\vec{T}_{q}, \vec{\bar{Q}}_{q}, \vec{R}_1 = \textsc{block-Lanczos}(\vec{A},\vec{\Omega},0,q)$
\State $\vec{Q} = \textsc{orth}(\vec{\bar{Q}}_{q}\left[f(\vec{T}_{q})\right]_{:,1:b}\vec{R}_1)$
\Comment{$\approx\textsc{orth}(f(\vec{A}) \vec{\Omega})$} \label{line:gsf:faOmega}
\State $\vec{T}_{n}, \vec{\bar{Q}}_{n} = \textsc{block-Lanczos}(\vec{A},\vec{Q},0,n)$
\State  $t_{\text{defl}} = \tr(\left[f(\vec{T}_{n})\right]_{1:b,1:b})$ 
\Comment{$\approx\tr(\vec{Q}^\T f(\vec{A})\vec{Q})$}
\State  $\vec{Y} = (\vec{I}-\vec{Q}\vec{Q}^\T)\vec{\Psi}$
\Comment{$\vec{Y} = [\vec{y}_1,\ldots,\vec{y}_{m}]$}
\For{$i=1,2,\ldots,m$}
\State $\vec{T}_{n}^{(i)} = \textsc{block-Lanczos}(\vec{A},\vec{y}_i,0,n)$ 
\State $t_{\text{rem}} = t_{\text{rem}} + \frac{d-b}{m} [f(\vec{T}_{n}^{(i)})]_{1,1}$
\Comment{$\frac{d-b}{m} {\vec{y}_i^\T \vec{B} \vec{y}_i}/{\vec{y}_i^\T \vec{y}_i}$} 
\EndFor
\State \Return $t_{\text{defl}} + t_{\text{rem}}$
\EndProcedure
\end{algorithmic}
\end{algorithm}

If $n$ and $q$ are large enough so that products with $f(\vec{A})$ are computed almost exactly, then approximation to $\tr(f(\vec{A}))$ output by \cref{alg:generalized_simple_func} should be close to that generated by \cref{alg:generalized_simple}.
The following result, which we prove in \cref{section:proofs}, guarantees that small errors in the computation of the projection matrix $\vec{Q}$ do not significantly impact the guarantees for \cref{alg:generalized_simple}.
The remaining impacts of the error in the approximation \cref{eqn:mv_approx} are analyzed formally in \cref{thm:main}.

\begin{lemma}\label{thm:hpp_func} 
Suppose \cref{alg:generalized_simple_func} terminates successfully.
Define 
\begin{align*}
    \Delta \equiv \| f(\vec{A}) \vec{\Omega} - \vec{\bar{Q}}_{q}\left[f(\vec{T}_{q})\right]_{:,1:b}\vec{R}_1 \|_2/\| \vec{\Omega} \|_2,
\end{align*}
where $\vec{\bar{Q}}_{q}$, $\vec{T}_q$, and $\vec{R}_1$ are as in \cref{alg:generalized_simple_func}.
Suppose $\vec{V} \equiv \textsc{orth}(f(\vec{A}) \vec{\Omega})$ and $\vec{Q} \equiv \textsc{orth}(\vec{\bar{Q}}_{q}\left[f(\vec{T}_{q})\right]_{:,1:b}\vec{R}_1)$ are full rank and introduce orthogonal projectors
\[
    \vec{P}_{\vec{V}} \equiv (\vec{I} - \vec{V} \vec{V}^\T)
    ,\qquad  \vec{P}_{\vec{Q}} \equiv (\vec{I} - \vec{Q} \vec{Q}^\T).
\]
Then,
\[
    \| \vec{P}_{\vec{Q}} f(\vec{A}) \vec{P}_{\vec{Q}} 
    -  \vec{P}_{\vec{V}} f(\vec{A}) \vec{P}_{\vec{V}} \|_{\F}
    \leq 2 \sqrt{2b} 
    \,\kappa(f(\vec{A}))\,\kappa(\vec{\Omega})
    \Delta.
\]
\end{lemma}

Note that while $\Delta$ is random, \cref{thm:func_rate} gives a deterministic prior bound for how large $q$ must be set, in terms of the best polynomial approximation to $f$ on $[\lambda_{\textup{min}},\lambda_{\textup{max}}]$, to ensure $\Delta$ is small (assuming the block-Lanczos algorithm \cref{alg:block_lanczos} terminates successfully).
We state the lemma in terms of $\Delta$ rather than the bound in \cref{thm:func_rate} as there are many bounds for Lanczos besides the simple ones stated in \cref{thm:func_rate}.

In almost all situations where sketching is used, the block size $b$ is much smaller than the dimension $d$. 
In such cases, the condition number of $\vec{\Omega}$ will be relatively small (e.g.~$<4$) with high probability \cite{vershynin2012introduction}.
Thus, \cref{thm:hpp_func} says that the projection step in \cref{alg:generalized_simple_func} behaves almost identically to that of \cref{alg:generalized_simple}, as long as $f(\vec{A}) \vec{\Omega}$ is computed fairly accurately in \cref{alg:generalized_simple_func}.

Note also that 
\[
    \| \vec{P}_{\vec{V}} f(\vec{A}) \vec{P}_{\vec{V}} \|_{\F}
    \leq \| f(\vec{A}) \vec{P}_{\vec{V}} \|_{\F}
    = \|f(\vec{A})  -  f(\vec{A}) \vec{V} \vec{V}^\T \|_{\F}.
\]
Bounds for the quality of the low-rank approximation $f(\vec{A}) \vec{V} \vec{V}^\T$ to $f(\vec{A})$ are well studied; see for instance \cite{meyer2021hutch,halko2011finding}.
In particular, $f(\vec{A}) \vec{V} \vec{V}^\T$ has nearly the same approximation error as the best rank-$k$ approximation to $f(\vec{A})$ in the Frobenius (or spectral) norm, at least when $b = k+p$ for some small $p$.

\section{Krylov-aware function approximation}
\label{sec:krylov_start}

The central insight of this paper is that instead of treating the product $\vec{B}\vec{\Omega}=f(\vec{A})\vec{\Omega}$ as a black-box routine, \cref{alg:generalized_simple_func} can be made more efficient by exploiting the structure of the block Krylov space built with $\vec{A}$ and $\vec{\Omega}$.
In particular, we make use of the fact that if $\vec{\bar{Q}}_{q+1}$ is an orthornomal basis for $\mathcal{K}_{q+1}(\vec{A},\vec{\Omega})$, then
\begin{equation}\label{eqn:krylov_aware}
    \mathcal{K}_n(\vec{A},\vec{\bar{Q}}_{q+1}) = \mathcal{K}_{q+n}(\vec{A},\vec{\Omega}).
\end{equation}
Indeed, we have that $\operatorname{span}(\vec{\bar{Q}}_{q+1}) 
    = \mathcal{K}_{q+1}(\vec{A},\vec{\Omega})
    = \operatorname{span}\{\vec{\Omega}, \vec{A}\vec{\Omega}, \ldots, \vec{A}^q \vec{\Omega} \}$.
Thus, 
\begin{align*}
    \mathcal{K}_n(\vec{A},\vec{\bar{Q}}_{q+1})
    &= \operatorname{span}\{\vec{\bar{Q}}_{q+1}, \vec{A}\vec{\bar{Q}}_{q+1}, \ldots, \vec{A}^{n-1} \vec{\bar{Q}}_{q+1} \}
    \\&= \operatorname{span}\{ \vec{\Omega}, \vec{A}\vec{\Omega}, \ldots, \vec{A}^q \vec{\Omega} , 
    \\&\hspace{5em}
    \vec{A} \vec{\Omega}, \vec{A}^2\vec{\Omega}, \ldots, \vec{A}^{q+1} \vec{\Omega}, 
    \\&\hspace{11em} \ddots
    \\
    & \hspace{7em} \vec{A}^q \vec{\Omega}, \vec{A}^{q+1}\vec{\Omega}, \ldots, \vec{A}^{q+n-1} \vec{\Omega}\}
    = \mathcal{K}_{n+q}(\vec{A},\vec{\Omega})
    .
\end{align*}
We note that the block-Lanczos algorithm \cref{alg:block_lanczos} can be implemented to automatically respect \cref{eqn:krylov_aware} without breaking down.\footnote{
Suppose that, when a rank-deficient matrix $\vec{Z}$ is encountered during the block-Lanczos process, the QR factorization in \hyperref[line:qnext]{line~\ref{line:qnext}} of \cref{alg:block_lanczos}} produces an orthonormal basis for the column space of $\vec{Z}$ and the block size is decreased. 
Then, in exact arithmetic, $q+n$ steps of block-Lanczos on the  orthonormal matrix $\vec{Q}_1 \in \mathbb{R}^{d\times b}$ produces the same output as $n$ steps of block-Lanczos on the input $\vec{\bar{Q}}_{q+1}\in \mathbb{R}^{d\times (q+1)b}$, where $\vec{\bar{Q}}_{q+1}$ is the result of running $q$ steps of block-Lanczos on $\vec{Q}_1$. 
We do not provide such an implementation as the details would obfuscate the main point of this paper.

The relation \cref{eqn:krylov_aware} suggests that, given $\vec{\bar{Q}}_{q+1}$, we can construct an approximation to $\vec{\bar{Q}}_{q+1}^\T f(\vec{A}) \vec{\bar{Q}}_{q+1}$ using $n-1$ additional products with matrices of size just $n\times b$.
In particular, we can use the approximation
\begin{equation}
    \vec{\bar{Q}}_{q+1}^\T \vec{\bar{Q}}_{q+n+1} f(\vec{T}_{q+n}) \vec{\bar{Q}}_{q+n+1}^\T\vec{\bar{Q}}_{q+1}
    = \left[f(\vec{T}_{q+n})\right]_{1:(q+1)b,1:(q+1)b}. \label{eqn:krylov_quad_approx}
\end{equation}
This approximation is exact when $f$ is a sufficiently low-degree polynomial. 
In particular, we have the following lemma.
\begin{lemma}\label{thm:krylov_quad_approx}
Assuming the Krylov subspace $\mathcal{K}_{q+n}$ has dimension  $(q+n)b$ so that no rank-deficient blocks are encountered, the approximation \eqref{eqn:krylov_quad_approx} is exactly equal to $\vec{\bar{Q}}_{q+1}^\T f(\vec{A}) \vec{\bar{Q}}_{q+1}$ whenever $f$ is a polynomial of degree at most $2n-1$.
\end{lemma}

In addition, the ``Krylov-aware'' approach yields a low-rank approximation to $f(\vec{A})$.
In particular, 
\begin{equation}
    \label{eqn:krylov_aware_low_rank}
    \vec{\bar{Q}}_{q+1} \left[f(\vec{T}_{q+n})\right]_{1:(q+1)b,1:(q+1)b} \vec{\bar{Q}}_{q+1}^\T
\end{equation}
gives an approximation to the symmetric projection of $f(\vec{A})$ onto $\mathcal{K}_{q+1}(\vec{A},\vec{\Omega})$ which is itself a low-rank approximation to $f(\vec{A})$.
A theoretical analysis of how well \cref{eqn:krylov_aware_low_rank} works for low-rank approximation would be an interesting topic of further study.
In particular, for functions $f$ such as the square root which compress the eigenvalues of $\vec{A}$, it seems reasonable that this approach would significantly outperform approaches based on approximating $f(\vec{A})\vec{\Omega}$ with a Krylov subspace method (see  \cite{persson2022randomized} for some related theoretical results in this direction).

\subsection{Main algorithm} \label{sec:greedy}

\begin{algorithm}[ht]
\caption{Krylov-aware stochastic trace estimation for matrix functions}\label{alg:main}
\fontsize{10}{14}\selectfont
\begin{algorithmic}[1]
\Procedure{krylov-aware-trace}{$\vec{A}, f,b,q, m, n$}
\State Sample Gaussian matrices $\vec{\Omega} \in \mathbb{R}^{d\times b}$ and $\vec{\Psi} \in \mathbb{R}^{d\times m}$
\State $\vec{T}_{q+n}, \vec{\bar{Q}}_{q+1} = \textsc{block-Lanczos}(\vec{A},\vec{\Omega},q,n)$ 
\label{line:greedy_blockLanczos}
\State $t_{\text{defl}} = \tr\left( \left[f(\vec{T}_{q+n})\right]_{1:(q+1)b,1:(q+1)b} \right)$\Comment{$\approx\tr(\vec{\bar{Q}}_{q+1}^\T f(\vec{A})\vec{\bar{Q}}_{q+1})$ }
\State $\vec{Y} = (\vec{I}-\vec{\bar{Q}}_{q+1}\vec{\bar{Q}}_{q+1}^\T)\vec{\Psi}$ \Comment{$\vec{Y} = [\vec{y}_1,\ldots,\vec{y}_{m}]$}
\For{$i=1,2,\ldots,m$} \label{line:greedy_forloop}
\State $\vec{T}_{n}^{(i)} = \textsc{block-Lanczos}(\vec{A},\vec{y}_i,0,n)$ \label{line:greedy_lanczos}
\State $t_{\text{rem}} = t_{\text{rem}} + \frac{d-(q+1)b}{m}[f(\vec{T}_{n}^{(i)})]_{1,1}$
\Comment{$\approx\frac{d-(b+1)q}{m} {\vec{y}_i^\T f(\vec{A}) \vec{y}_i}/{\vec{y}_i^\T \vec{y}_i}$}  \label{line:greedy_endfor}
\EndFor
\State \Return $t_{\text{defl}} + t_{\text{rem}}$
\EndProcedure
\end{algorithmic}
\end{algorithm}
Our ``Krylov-aware'' approach to stochastic trace estimation is summarized in \cref{alg:main}. 
We make a few observations about the implementation: 
\begin{itemize}[itemsep=0pt]
    \item In \hyperref[line:greedy_blockLanczos]{line \ref{line:greedy_blockLanczos}}, the columns of $\vec{\bar{Q}}_{q+1}$ should be kept orthonormal. The remaining columns of $\vec{\bar{Q}}_{q+n+1}$ do not need to be reorthogonalized against $\vec{\bar{Q}}_{q+1}$.
    \item In \hyperref[line:greedy_lanczos]{line \ref{line:greedy_lanczos}}, the Krylov basis vectors do not need to be reorthogonalized.
    \item Lines \ref{line:greedy_forloop}--\ref{line:greedy_endfor} can be run in parallel, blocking matvecs with $\vec{A}$. While it would also be possible to run $\textsc{block-Lanczos}(\vec{A},\vec{Y},0,n)$, doing so does not seem to significantly improve the quality of the estimate since the columns of $\vec{Y}$ do not share any notable structure with respect to $\vec{A}$.
    \item Lines \ref{line:greedy_forloop}--\ref{line:greedy_endfor} can also be run in parallel with the final $n$ iterations of \hyperref[line:greedy_blockLanczos]{line \ref{line:greedy_blockLanczos}}. 
    \item The algorithm can easily be adapted to approximate $\tr(f(\vec{A}))$ for many functions $f$ (which do not even need to be known in advance) at minimal additional cost.
\end{itemize}

We can provide a simple error guarantee for \cref{alg:generalized_simple_func,alg:main}, given a prescribed choice of $b$, $q$, and $n$.
We discuss heuristics and intuition for how $b$ and $q$ can be chosen in the next section.
\begin{theorem}\label{thm:main}
Assume that $\vec{A}$, $f$, $b$, $q$, $n$ are such that \cref{alg:generalized_simple_func,alg:main} terminate successfully with probability one.
Set $\vec{\hat{Q}} \equiv \vec{Q}$ in the case of \cref{alg:generalized_simple_func} and $\vec{\hat{Q}} \equiv \vec{\bar{Q}}_{q+1}$ in the case of \cref{alg:main}, and define $\hat{b}$ as the number of columns in $\vec{\hat{Q}}$.
Define the random variable $\Delta$ by
\begin{align*}
    \Delta \equiv \max\{ \Delta_0, \Delta_1,  \ldots, \Delta_m\}
    ,\qquad
    \begin{cases}
    \Delta_0&\equiv \big\| \vec{\hat{Q}}^\T f(\vec{A}) \vec{\hat{Q}} - \left[f(\vec{{T}}_{q+n})\right]_{1:\hat{b},1:\hat{b}} \big\|_2\big/\| \vec{\hat{Q}} \|_2^2
\\
    \Delta_i& \equiv \big| {\vec{y}_i^\T f(\vec{A})\vec{y}_i} - \|\vec{y}_i \|_2^2[f(\vec{{T}}_{n}^{(i)})]_{1,1} \big| \big/ \|\vec{y}_i \|_2^2\end{cases}.
\end{align*}
Here $\vec{T}_{q+n}$, $\vec{T}_n^{(i)}$, and $\vec{y}_i = (\vec{I} -\vec{\hat{Q}}\vec{\hat{Q}}^\T) \bm{\uppsi}_i$ are as in each algorithm.
Then, the outputs of the \cref{alg:generalized_simple_func,alg:main} satisfy
\[
    \left| \tr(f(\vec{A})) - \EE\big[ (t_{\textup{defl}} +t_{\textup{rem}})  \big]\right| \leq d \, \EE\big[\Delta \big],
\]
\[
    \VV\big[ t_{\textup{defl}} +t_{\textup{rem}} \big]
    \leq \Big(\sqrt{V} + d\sqrt{\EE\big[\Delta^2  \big]}\Big)^2,
\]
where, 
with $\vec{F} \equiv (\vec{I} - \vec{\hat{Q}}\vec{\hat{Q}}^\T) f(\vec{A}) (\vec{I} - \vec{\hat{Q}}\vec{\hat{Q}}^\T)$, 
\[
    V \equiv \frac{2}{m} \left(\frac{d-\hat{b}}{d-\hat{b}+2}\right)\left(  \EE\big[\|\vec{F}\|_\F^2 \big] - \frac{\EE\big[\tr(\vec{F})^2 \big]}{d-\hat{b}} \right).
\]
Moreover, with $\vec{P}_{\hat{\vec{Q}}} \equiv \vec{I} - \hat{\vec{Q}}\hat{\vec{Q}}^\T$, and provided the same matrix $\vec{\Omega}$ is used in both algorithms, 
\( \| \vec{P}_{\vec{\bar{Q}}_{q+1}} f(\vec{A}) \vec{P}_{\vec{\bar{Q}}_{q+1}} \|_\F^2
\leq \| \vec{P}_{\vec{Q}} f(\vec{A}) \vec{P}_{\vec{Q}} \|_\F^2 \) 
 and so 
\[
\EE\big[\| \vec{P}_{\vec{\bar{Q}}_{q+1}} f(\vec{A}) \vec{P}_{\vec{\bar{Q}}_{q+1}} \|_\F^2 \big]
\leq \EE\big[\| \vec{P}_{\vec{Q}} f(\vec{A}) \vec{P}_{\vec{Q}} \|_\F^2  \big] .\]
\end{theorem}

We make several comments on \cref{thm:main}.
First assuming the algorithm terminates successfully, $\Delta$ can be upper bounded by a deterministic quantity which, for reasonable functions $f$, tends to zero as $n\to\infty$; see \cref{bgrd:krylov}.
Second, the variance bound is obtained using the fact that $\frac{d-\hat{b}}{m} {\vec{y}_i^\T f(\vec{A}) \vec{y}_i}/{\vec{y}_i^\T \vec{y}_i}$ has known variance $V$ (with respect to the randomness in $\bm{\uppsi}_i$)\cite{girard1987algorithme,martinsson2020randomized}. 
The concentration of measure phenomenon for the uniform distribution on the sphere \cite{ledoux2001concentration} implies that this random variable is sub-Gaussian \cite{popescu2006entanglement}. 
However, fine-grained concentration inequalities are are not as readily available in the literature as the inequalities for Gaussian or Rademacher vectors \cite{meyer2021hutch,persson2022improved}.
Finally, the result shows that if $d\Delta$ is small relative to $\max(1,
\sqrt{V})$, then \cref{alg:main} cannot perform significantly worse than \cref{alg:generalized_simple_func}.
Our experiments in \cref{sec:experiments} show that \cref{alg:main} often performs significantly better than \cref{alg:generalized_simple_func}.

\Cref{thm:hpp_func} in conjunction with \cref{thm:func_rate} implies that when $q$ increases, the quality of the projection matrix $\vec{Q}$ by \cref{alg:generalized_simple_func} converges to the projection matrix which would be generated if products with $f(\vec{A})$ were computed exactly (i.e.~to the projection matrix generated by \hpp). 
For many cases of interest, the latter is a poor way to spend one's computational budget. For example, if $\vec{A}$ is symmetric positive semidefinite and $f(x) = x^{1/2}$, then a better projection matrix $\vec{Q}$ could be obtained simply by running one step of subspace iteration with $\vec{A}$ itself! 
On the other hand, one can expect the block Krylov subspace $\vec{\bar{Q}}_{q+1}$ to continue to improve as $q$ increases.\footnote{If $\vec{A}$ has only a few distinct eigenvalues, the dimension of the block-Krylov subspace may stop growing. However, in such cases the approximation $f(\vec{A})\vec{\Omega}$ is exact, so the projection space used by \cref{alg:main} is no worse than what is used by \cref{alg:generalized_simple_func} (which also becomes exact and therefore is the same as the space used by \cref{alg:generalized_simple}.)} Exactly how much the block Krylov space improves, however, is a more difficult question to answer. We discuss the topic further in \cref{sec:past_work}.

\subsection{Choice of parameters}\label{sec:parameter_selection}
In all, \cref{alg:main} requires $b(q+n) + mn$ matvecs with $\vec{A}$ and a minimum of $q+n$ matrix loads, matching the costs of \cref{alg:generalized_simple_func}.

This interpretation alone is perhaps too optimistic, since \cref{alg:main} requires storing and orthogonalizing $\vec{\bar{Q}}_{q+1}$. 
In many situations, the size of $\vec{\bar{Q}}_{q+1}$ will be limited by these storage and reorthogonalization costs rather than the number of matvecs with $\vec{A}$. 
Suppose that we have some predetermined limit on the number of vectors to be used for deflation, i.e.~$b(q+1)$ is bounded by a constant. What block size will be most useful? Using a small block size will reduce the number of matvecs, while using a larger block size will allow for greater parallelism and reduce the number of matrix loads.

\begin{table}[h]
    \caption{Approximate costs of algorithms. Here storage is measured in terms of the number of vectors of length $d$ stored.
    For simplicity, we assume each implicit product with $f(\vec{A})$ in \hpp{} is sequentially implemented using $n$ steps of the Lanczos method for matrix function approximation.\protect\footnotemark 
    The number of matrix loads of $\vec{A}$, which is more important than the total number of matvecs in some situations, is not displayed here. In general, the number of matrix loads can be reduced at the cost of increased storage.
    }
    \label{tab:approx_costs}
    \centering
    \begin{tabular}{ccccc}\toprule
    Algorithm & number matvecs & storage & inner products \\ \midrule
    \hyperref[alg:generalized_simple_func]{Alg.~\ref{alg:generalized_simple_func}} &
   $b(q+n) + mn$ 
    &  $O(b+n)$
    & $O(b^2+bn+bm+mn)$ 
    \\
    \hyperref[alg:main]{Alg.~\ref{alg:main}} & $b(q+n)+mn$ & $O(bq+n)$ & $O(b^2(q^2+n)+bqm+mn)$ 
    \\ 
    \hyperref[alg:restart]{Alg.~\ref{alg:restart}} & $b(qr+q+n)+mn$ & $O(bq+n)$ & $O(b^2(q^2r+n)+bqm+mn)$\\ 
    \bottomrule
    \end{tabular}

\end{table}
\footnotetext{It is possible to reduce the storage costs of Lanczos-based matrix function approximation from $O(n)$ to $O(1)$ using a two-pass approach. However, this doubles the number of matrix-vector products from $n$ to $2n$. Other low memory approaches may be possible for certain functions \cite{guttel2021comparison}.}

Methods focused on computing a small number of eigenvectors to high accuracy have tended to use small block sizes (including $b=1$) \cite{wu2000thick,stewart2002krylov} while more recent analyses of randomized block-Lanczos have tended to consider large $b$ and small $q$ \cite{musco2015randomized,wang2015improved,drineas2018structural,halko2011algorithm}; see \cite{martinsson2020randomized} for more background. 
For the problems of finding a low-rank approximation or a few eigenpairs of $\vec{A}$, the total number of matrix loads is roughly proportional to $q$, so the trade-off between the Krylov-depth $q$ and block-size $b$ is not just in the total number of matvecs, but also the number of matrix-loads.

In our situation the minimum number of matrix loads is $q+n$.
Since $n$ will typically be relatively high (i.e.~not some small constant like 3) due to the need to construct a sufficiently large Krylov subspace such that products with $f(\vec{A})$ are computed to sufficient accuracy, making $q$ small cannot significantly reduce the overall number of matrix-loads.
In other words, the presumed benefit of our algorithm is in the reduced number of matvecs required, rather than a reduction in the number of matrix loads.
This suggests that when the size of the deflation space is limited, it makes sense to take $b$ small and $q$ large.

To improve the quality of the approximation to the dominant eigenspace of $f(\vec{A})$ while limiting the size of $\vec{\bar{Q}}_{q+1}$ (and therefore reorthogonalization and storage costs) it is common to use restarting schemes for block-Lanczos \cite{baglama2003implicitly,zhou2008block}.
We describe how these approaches can be used with our algorithm in \cref{sec:low_mem}. The costs of the algorithms described in this paper are summarized in \cref{tab:approx_costs}. One feature worth noting is that since we assume the last $n$ steps of the block-Lanczos process are done without reorthogonalization, the number of inner products required by our algorithms grow at worst proportional to $n$ (as opposed to $n^2$).

\subsection{Relation to past work}
\label{sec:past_work}

A number of papers have integrated Krylov subspace methods into implicit trace approximation algorithms.
The most widespread approach has been to use Krylov subspace methods to approximate products with $f(\vec{A})$ for quadratic trace estimators \cite[etc.]{bai1996largescale,bai1996bounds,skilling1989eigenvalues,weisse2006kernel,han2017approximating,ubaru2017fast}; see \cite{chen2022randomized} for a recent review.
This approach is well understood and has been studied in detail.
More recently, a number of algorithms aim to incorporate low-rank approximation.

The paper \cite{li2021randomized} takes $\vec{Q}$ to be an orthonormal basis for the entire block Krylov subspace $\mathcal{K}_{q+1}(\vec{A},\vec{\Omega})$ and uses it to estimate $\tr(\vec{A})$ and $\tr(\log(\vec{I} + \vec{A}))$ when $\vec{A}$ is positive definite. This method builds on the work in \cite{saibaba2017randomized} (which uses randomized subspace iteration) and is a special case of \cref{alg:main} with $n=1$ and $m=0$. 
Our approach works with any $n$ which is critical for general matrix functions whose dominant eigenvalues may not align with those of $\vec{A}$. 

The paper \cite{wu2016estimating} estimates $\tr(\vec{A}^{-1})$ using a low-rank approximation to $\vec{A}$ for variance reduction. The authors note that the Krylov space produced with each new sample of the quadratic trace estimator contains information about the eigenvalues and eigenvectors of $\vec{A}$, and discuss the possibility of using this information to incrementally update the low-rank approximation. 
Low-rank approximation and quadratic trace estimation were combined in \cite{morita2020finite} in order to approximate $\tr(\exp(-\beta \vec{A}))$ and related quantities.

A more direct inspiration for our work is \cite{bekas2007estimator}, in which the authors propose to estimate the diagonal of a matrix using as samples the vectors $\{ \vec{v}_0,\ldots,\vec{v}_k \}$, where $\vec{v}_k \equiv T_k(\vec{A})\vec{v}_0$ and $T_k$ is the $k$-th Chebyshev polynomial of the first kind. By exploiting the shared structure of these vectors, they significantly reduce computational costs. The authors note that these samples are ``correlated and not completely random''; our innovation is to apply their technique to the deflation step only, which will allow the trace estimator to remain unbiased. 

All of \cite{saibaba2017randomized,li2021randomized,persson2022randomized} derive a priori bounds for low-rank approximation of matrix functions. 
In particular, \cite{persson2022randomized} analyzes low rank approximation of operator monotone functions using a Nystr\"om based approximation similar to what is used by a subcase of our restarted variant described in \cref{sec:low_mem} (again with $n=1$).
This approach is combined with quadratic trace estimation to produce estimates of the trace of operator monotone functions.
All of these analyses make critical use of the structure of $f$, which allows $\vec{\hat{Q}}\vec{\hat{Q}}^\T f(\vec{A}) \vec{\hat{Q}}\vec{\hat{Q}}^\T$ to be approximated by $f(\vec{\hat{Q}} \vec{\hat{Q}}^\T \vec{A} \vec{\hat{Q}} \vec{\hat{Q}}^\T)$ if $\vec{A}$ is well-approximated by $\vec{\hat{Q}} \vec{\hat{Q}}^\T \vec{A} \vec{\hat{Q}} \vec{\hat{Q}}^\T$.
For arbitrary $f$, it seems possible that $\vec{\hat{Q}}\vec{\hat{Q}}^\T f(\vec{A}) \vec{\hat{Q}}\vec{\hat{Q}}^\T$ provides a good approximation to $f(\vec{A})$ even if $\vec{\hat{Q}} \vec{\hat{Q}}^\T \vec{A} \vec{\hat{Q}} \vec{\hat{Q}}^\T$ does not provide a good approximation to $\vec{A}$.
Thus, the observation that $\vec{\hat{Q}}^\T f(\vec{A}) \vec{\hat{Q}}$ can be approximated efficiently if the span of the columns of $\vec{\hat{Q}}$ is itself a Krylov subspace (see \cref{eqn:krylov_aware_low_rank}) is a critical aspect of our algorithm.
However, this added complexity means that a priori bounds for the quality of the ``Krylov-aware'' low-rank approximation used in our algorithm are outside of the scope of this paper.
Understanding the theoretical behavior of these algorithms is an interesting topic for further study

\section{Variants} \label{sec:variants}
We discuss several variants of \cref{alg:main} which may be more practical in many situations.

\subsection{Adaptive variant}\label{sec:adaptive}

One disadvantage of \cref{alg:main} is that a user may not know ahead of time how much effort should be devoted to variance reduction as opposed to the quadratic trace estimator. In order to ameliorate this problem, we propose an adaptive variant (\cref{alg:adaptive}) that takes parameters $\epsilon > 0$ and $\delta\in (0,1)$ and attempts to compute an estimate satisfying 
\[
 \mathbb{P}\left( |\mathsf{est} - \tr(f(\vec{A}))| > \epsilon \right) < \delta.
\]
The adaptive algorithm \cite[Algorithm 2]{persson2022improved} provably satisfies such an error guarantee under the assumption that matvecs with $f(\vec{A})$ can be computed exactly. 
We model our algorithm instead on the slightly less rigorous A-\hpp{} \cite[Algorithm 3]{persson2022improved}, which is simpler while still having good empirical behavior. 

We assume that the block-Lanczos process computes matvecs with $f(\vec{A})$ exactly\textemdash or more precisely, we assume that $n$ is set large enough for the discrepancy to be ignored. For more detailed discussions on how to select the degree $n$ depending on the function $f$ and spectrum of $\vec{A}$, see \cite{golub2009matrices,ubaru2017fast,chen2022randomized}.

We also assume that the block size $b$ is fixed. Once this is done, the adaptive algorithm must decide first how many Lanczos iterations $q$ to devote to variance reduction, and second how many vectors $m$ to use for the quadratic trace estimator. The answer to the second part is straightforward, given an estimate of the remainder $\|\widetilde{\vec{B}}\|_\F$. 
In \cite{cortinovis2021randomized} it is shown that for user-specified tolerances $\epsilon>0$ and $\delta \in (0,1)$, 
\begin{equation}
    \PP\left[ \big| \bm{\uppsi}^\T \widetilde{\vec{B}} \bm{\uppsi} -\tr(\widetilde{\vec{B}})\big| > \epsilon \right] < \delta 
    \quad
    \text{if}
    \quad m \geq c\,\epsilon^{-2}(\|\widetilde{\vec{B}}\|_\F^2 + \epsilon\|\widetilde{\vec{B}}\|_2)\ln(2/\delta),\label{eqn:epsDelta}
\end{equation}
where $c$ is a constant depending on the distribution of $\bm{\uppsi} \in \mathbb{R}^{d\times m}$.\footnote{It suffices to take $c=4$ for Gaussian vectors  and $c=16$ for Rademacher \cite{cortinovis2021randomized}.}

Relaxing the bound \eqref{eqn:epsDelta} via the approximation $\|\widetilde{\vec{B}}\|_{\F}^2 + \epsilon \|\widetilde{\vec{B}}\|_2 \approx \|\widetilde{\vec{B}}\|_{\F}^2$, the algorithm A-\hpp{} defines 
\begin{equation}\label{def:ced}
 C(\epsilon,\delta) \equiv 4\epsilon^{-2}\log(2/\delta)
\end{equation}
and proposes to use $m = C(\epsilon,\delta)\|\widetilde{\vec{B}}\|_\F^2$ samples. In practice A-\hpp{} increases $m$ incrementally while simultaneously updating an estimate of $\|\widetilde{\vec{B}}\|_\F^2$, and we do the same in \cref{alg:adaptive}. We use a slightly paraphrased version of \cite[Lemma 2.2]{persson2022improved}, which relies on a result from \cite{roosta2015assessing}.
\begin{lemma}\label{lemma:alpha}
 Let $\vec{\Psi}\in \mathbb{R}^{d\times k}$ be a standard Gaussian matrix, and let $\widetilde{\vec{B}}\in \mathbb{R}^{d\times d}$. For any $\alpha \in (0,1)$, it holds that
 \begin{equation}
     \mathbb{P}\left(\frac{1}{k\alpha}\|\widetilde{\vec{B}}\vec{\Psi}\|_\F^2\leq \|\widetilde{\vec{B}}\|_\F^2\right) \leq \mathbb{P}(X\leq k\alpha), 
 \end{equation}
 where $X\sim \chi_k^2$ is a chi-squared random variable with $k$ degrees of freedom. 
\end{lemma}
For a user-specified failure tolerance $\delta$, we can therefore define 
\begin{equation}\label{def:alpha}
    \alpha_k \equiv \frac{1}{k}F_X^{-1}(\delta),
\end{equation}
where $F_X$ is the cumulative distribution function (CDF) of $X\sim\chi_k^2$. It follows from \cref{lemma:alpha} that $\tfrac{1}{k\alpha_k}\|\widetilde{\vec{B}}\vec{\Psi}\|_\F^2$ will overestimate $\|\widetilde{\vec{B}}\|_\F^2$ with probability at least $1-\delta$. It is shown in \cite{persson2022improved} that the sequence $\{\alpha_k\}$ increases monotonically and converges to 1. \cref{alg:adaptive} uses the same set of random vectors to incrementally estimate $\tr(\vec{\widetilde{B}})$ and $\|\vec{\widetilde{B}}\|_\F^2$ simultaneously, stopping when the number of samples $k$ exceeds the sampling number
\[
    m_k \equiv \frac{1}{k\alpha_k}\|\vec{\widetilde{B}}\vec{\Psi}\|_\F^2 \gtrapprox \|\vec{\widetilde{B}}\|_\F^2. 
\]

\begin{remark}
 The factor $\alpha_k$ is pessimistic if the stable rank  $\|\vec{\widetilde{B}}\|_\F^2/\|\vec{\widetilde{B}}\|_2^2$ is large. Ideally, one should incrementally estimate the stable rank and use it in combination with bounds such as those in \cite{gratton2018improved,cortinovis2021randomized} to obtain values for $\alpha_k$ closer to 1, and therefore smaller values for $m_k$. We stick with the definition \eqref{def:alpha} in order to focus on the more novel aspects of our algorithm. 
\end{remark}

All that remains is to determine the number of iterations $q$. The total number of matvecs used by \cref{alg:main}, which we use as a proxy for the computational cost, is 
\begin{align*}
    M(q) &\equiv (q+n)b + nm \\
        &\approx (q+n)b + nC(\epsilon,\delta)\|(\vec{I}-\vec{\bar{Q}}_{q+1}\vec{\bar{Q}}_{q+1}^\T)f(\vec{A})(\vec{I}-\vec{\bar{Q}}_{q+1}\vec{\bar{Q}}_{q+1}^\T)\|_\F^2\\
        &= (q+n)b + nC(\epsilon,\delta)\big(\|f(\vec{A})\|_\F^2 - 2\|f(\vec{A})\vec{\bar{Q}}_{q+1}\|_\F^2 + \|\vec{\bar{Q}}_{q+1}^\T f(\vec{A})\vec{\bar{Q}}_{q+1}\|_\F^2\big)\\
        &= qb - nC(\epsilon,\delta)\big(2\|f(\vec{A})\vec{\bar{Q}}_{q+1}\|_\F^2 - \|\vec{\bar{Q}}_{q+1}^\T f(\vec{A})\vec{\bar{Q}}_{q+1}\|_\F^2\big) + \text{constant}. 
\end{align*}
Using $n$ steps of block-Lanczos to estimate the above quantity, we get the objective function
\[
    \widetilde{M}(q)\equiv qb - nC(\epsilon,\delta)\big(2\|[f(\vec{T}_{q+n})]_{:,1:(q+1)b}\|_\F^2-\|[f(\vec{T}_{q+n})]_{1:(q+1)b,1:(q+1)b}\|_\F^2\big).
\]
As is done in \cite[\S 2.1.2]{persson2022improved}, we propose to stop when $\widetilde{M}(q) > \widetilde{M}(q-1)> \widetilde{M}(q-2)$, which suggests that a local minimum has been found.\footnote{In practice, we may also include a preset limit $q_{\text{max}}$ on the number of vectors we are willing to store and orthogonalize. In such cases, one might turn to restarting (as described in the next section) to refine the quality of the low-rank approximation.}
With this approach the quantity $\widetilde{M}(q)$ can only be estimated after a delay of $n$ further iterations, but since \cref{alg:main} already requires the computation of $\vec{T}_{q+n}$ this work is not wasted. Our approach also requires us to compute the eigenvalues of the block tridiagonal matrix $\vec{T}$ at each iteration, but in theory the eigenvalues can be updated cheaply with each iteration.

In practice, the true cost of the algorithm is not necessarily directly proportional to the number of matvecs. 
For instance, storage or reorthogonalization costs are often a computational bottleneck.
In such cases, one can update the function $M(q)$ to account for such costs. 
However, this may make the resulting optimization problem involving $\widetilde{M}(q)$ more difficult.

\begin{algorithm}[ht]
\caption{Adaptive Krylov-aware stochastic trace estimation}\label{alg:adaptive}
\fontsize{10}{14}\selectfont
\begin{algorithmic}[1]
\Procedure{ada-trace}{$\vec{A}, f, n, b, \epsilon, \delta$}
\State Sample Gaussian $\vec{\Omega}\in \mathbb{R}^{d\times b}$
\State  $\vec{T}_{q+n}, \vec{\bar{Q}}_{q+1} = \textsc{block-Lanczos}(\vec{A},\vec{\Omega})$, run until $\widetilde{M}(q)$ minimized
\State $t_{\text{defl}} = \tr\left( \left[f(\vec{T}_{q+n})\right]_{1:(q+1)b,1:(q+1)b}\right)$
\State Initialize $t_{\text{rem}} = 0$ and $t_{\text{fro}} = 0$ \label{line:ada_remStart}
\State Initialize $m_0 = \infty$ and $k = 0$
\While {$m_k > k$}
\State $k = k+1$
\State Sample Gaussian $\boldsymbol{\uppsi}_k\in\mathbb{R}^d$
\State Orthogonalize $\vec{y}_k = (\vec{I} - \vec{\bar{Q}}_{q+1}\vec{\bar{Q}}_{q+1}^\T)\boldsymbol{\uppsi}_k$
\State $\vec{T}_n^{(k)} = \textsc{block-Lanczos}(\vec{A},\vec{y}_k,n)$
\State Increment $t_{\text{rem}} = t_{\text{rem}} + [f(\vec{T}_n^{(k)})]_{1,1}\|\vec{y}_k\|_2^2$ \Comment{Estimate $\tr(\vec{\widetilde{B}})$}
\State Increment $t_{\text{fro}} = t_{\text{fro}} + \|[f(\vec{T}_n^{(k)})]_{:,1}\|_2^2\,\|\vec{y}_k\|_2^2$ \Comment{Estimate $\|\vec{\widetilde{B}}\|_\F^2$}\label{line:ada_fnorm}
\State $\alpha_k = \frac{1}{k}F_X^{-1}(\delta)$, where $X\sim \chi_k^2$\Comment{Inverse CDF}\label{line:ada_gammafun}
\State $m_k = \tfrac{1}{k\alpha_k}C(\epsilon,\delta)t_{\text{fro}}$ \label{line:ada_sampling}
\EndWhile
\Return $t_{\text{defl}} + \frac{1}{k}t_{\text{rem}}$
\EndProcedure
\end{algorithmic}
\end{algorithm}

\cref{alg:adaptive} presents the adaptive algorithm. As with \cref{alg:main}, the loop can be blocked for efficiency. 
Note here that the normalization of $\vec{y}_k$ differs from \cref{alg:main}. 
This is essentially an artifact of analysis, as fine-grained concentration inequalities depending on $\|\widetilde{\vec{B}}\|_\F$ for vectors sampled from the hypersphere are less readily available. 
It is trivial and inexpensive to maintain both normalizations.

\subsection{Low-memory variant}
\label{sec:low_mem}

In practice, the memory and orthogonalization costs of block-Lanczos may limit the dimension of the Krylov subspace used for variance reduction. 
In such situations, we aim to find a  subspace of $\mathcal{K}_{q+1}(\vec{A},\vec{\Omega})$ that approximates the dominant eigenspace of $f(\vec{A})$ without having to store a basis for the entire block Krylov space. Some of the most widely used techniques for this task are based on restarting the Lanczos recurrence using a carefully chosen subspace of $\mathcal{K}_{q}(\vec{A},\vec{\Omega})$ \cite{baglama2003implicitly,sorensen1992implicit,zhou2008block}.
For concreteness and clarity we focus on the implicitly restarted block-Lanczos method \cite{baglama2003implicitly}.
Other techniques may make more sense in specific situations.

\begin{algorithm}[ht]
\caption{Low-memory Krylov-aware stochastic trace estimation}\label{alg:restart}
\fontsize{10}{14}\selectfont
\begin{algorithmic}[1]
\Procedure{restart-trace}{$\vec{A},f,r,\{p^{(i)}\},b,q, m, n$}
\State Sample Gaussian matrices $\vec{\Omega} \in \mathbb{R}^{d\times b}$ and $\vec{\Psi} \in \mathbb{R}^{d\times m}$
\For{$i=1,2,\ldots,r$}
\State $\vec{T}_{q}, \vec{\bar{Q}}_{q+1} = \textsc{block-Lanczos}(\vec{A},\vec{\Omega},q)$  \label{line:restart1}
\State $\vec{\Omega} = \vec{\bar{Q}}_q [ p(\vec{T}_q)]_{:,1:b} \vec{R}_1$ \label{line:restart2} \Comment{$\vec{\Omega} = p^{(i)}(\vec{A}) \vec{\Omega}$}
\EndFor
\State $\vec{T}_{q+n}, \vec{\bar{Q}}_{q+1} = \textsc{block-Lanczos}(\vec{A},\vec{\Omega},q+n)$ 
\State $t_{\text{defl}} = \tr\left( \left[f(\vec{T}_{q+n})\right]_{1:(q+1)b,1:(q+1)b} \right)$\Comment{$\tr(\vec{\bar{Q}}_{q+1}^\T f(\vec{A})\vec{\bar{Q}}_{q+1})$} \label{line:restartdefl}
\State $\vec{Y} = (\vec{I}-\vec{\bar{Q}}_{q+1}\vec{\bar{Q}}_{q+1}^\T)\vec{\Psi}$ \Comment{$\vec{Y} = [\vec{y}_1,\ldots,\vec{y}_{m}]$}
\For{$i=1,2,\ldots,m$} 
\State $\vec{T}_{n}^{(i)} = \textsc{block-Lanczos}(\vec{A},\vec{y}_i,n)$ 
\State $t_{\text{rem}} = t_{\text{rem}} + \frac{d-(q+1)b}{m}[f(\vec{T}_{n}^{(i)})]_{1,1}$
\Comment{$\frac{d-(q+1)b}{m} {\vec{y}_i^\T \vec{B} \vec{y}_i}/{\vec{y}_i^\T \vec{y}_i}$}
\EndFor
\State \Return $t_{\text{defl}} + t_{\text{rem}}$
\EndProcedure
\end{algorithmic}
\end{algorithm}

Specifically, suppose the block-Lanczos process has been run for $q$ iterations to obtain $\{\vec{Q}_k\}, \{\vec{M}_k\}, \{ \vec{R}_k \}$ which satisfy
\[
    \vec{A} \vec{\bar{Q}}_q = \vec{\bar{Q}}_q \vec{T}_q + \vec{Q}_{q+1} \vec{R}_{q+1} \vec{E}_q^\T.
\]
The restarting process involves updating $\vec{\Omega} = p(\vec{A}) \vec{\Omega} = \vec{\bar{Q}}_q [ p(\vec{T}_q)]_{:,1:b} \vec{R}_1$, where $p$ is some degree-$q$ polynomial, and then generating a new Lanczos recurrence starting with the updated $\vec{\Omega}$.
This process can then be repeated.
Ideally $p$ is large on the desirable eigenvalues of $\vec{A}$ and small elsewhere; i.e.~$p$ acts as a filtering polynomial.
A number of techniques for choosing such polynomials and performing the update step have been studied \cite{baglama2003implicitly,sorensen1992implicit}.

If filter polynomials $\{ p^{(i)} \}_{i=1}^{r}$ are used, then the final Krylov subspace generated is
\[
    \mathcal{K}_{q+1}(\vec{A}, p^{(r)}(\vec{A}) \cdots p^{(1)}(\vec{A}) \vec{\Omega}).
\]
While $p^{(r-1)}(\vec{A}) \cdots p^{(1)}(\vec{A}) \vec{\Omega}$ could be computed explicitly, the advantage of breaking it into phases is that the process can adapt to the information gained at each step; i.e.~the number of restarts $r$ and the polynomials $\{p^{(i)}\}_{i=1}^r$ do not have to be chosen in advance. 
This allows the process to be terminated after a sufficiently desirable approximation is obtained.

One might hope to apply the techniques from \cref{sec:adaptive} in order to derive a stopping criterion for the restarting stage of \cref{alg:restart}.
The main difference is that in many of the use cases we envision for \cref{alg:restart}, the rank of the approximation space $\vec{\bar{Q}}_{q+1}$ will be fixed.
In such cases, the minimum possible value of $\|\widetilde{\vec{B}}\|_\F^2 = \|(\vec{I}-\vec{\bar{Q}}_{q+1}\vec{\bar{Q}}_{q+1}^\T)f(\vec{A})(\vec{I}-\vec{\bar{Q}}_{q+1}\vec{\bar{Q}}_{q+1}^\T)\|_\F^2$ will be limited by the quality of the best rank $(q+1)b$ approximation to $f(\vec{A})$.
Thus, it would be more appropriate to terminate the restarting procedure after the quality of approximation is not improved.
This might be done by observing the quantity
\[
    \tr\left( \left[f(\vec{T}_{q+1})\right]_{1:(q+1)b,1:(q+1)b} \right)
\]
before deciding whether to repeat \hyperref[line:restart1]{lines \ref{line:restart1}} and \ref{line:restart2} or proceed to \hyperref[line:restartdefl]{line \ref{line:restartdefl}}.

One can of course use \cref{lemma:alpha} to estimate $\|\tilde{\vec{B}}\|_\F^2$ in order to determine the number of samples to be used in the second stage of the algorithm.
However, the spread of the output of the second stage of the algorithm can be easily controlled by relatively simple statistical methods for \emph{scalar} random variables such as observing the sample variance, bootstrapping, or jackknife.

While we believe that a practical version of our algorithm for application to real problems should incorporate both adaptive parameter selection and restarting, a good implementation is necessarily dependent on the problem at hand and the computing system to be used. Thus, for the sake of clarity, we do not provide a description which combines both  adaptive parameter selection and restarting.

\section{Numerical experiments} \label{sec:experiments}

\subsection{Quantum spin systems}

In this example, we consider the task of computing the partition function 
\begin{align*}
    Z(\beta) = \tr(\exp(-\beta \vec{A}))
\end{align*}
for the isotropic XY Heisenberg spin chain with a magnetic field of strength $h$ pointed in the $\textup{z}$-direction \cite{weisse2006kernel,schnack2020accuracy,faridfar2021thermodynamic,chen2022randomized}.
Specifically, the Heisenberg Hamiltonian for such a chain with \( N \) spins of spin number \(s=1/2\) is given by 
\begin{align*}
    \vec{A} = 
    2\sum_{i=1}^{N-1} \left(  \vec{s}^{\textup{x}}_i  \vec{s}^{\textup{x}}_{i+1} 
    +  \vec{s}^{\textup{y}}_i  \vec{s}^{\textup{y}}_{i+1} \right)
    + h \sum_{i=1}^{N} \vec{s}^{\textup{z}}_i.
\end{align*}
Here \( \vec{s}^{\textup{x}/\textup{y}/\textup{z}}_i  \in \mathbb{C}^{(2s+1)^N\times (2s+1)^N} \) is defined by
\begin{align*}
    \vec{s}^{\textup{x}/\textup{y}/\textup{z}}_i
    = \underbrace{\vec{I} \otimes \cdots \otimes \vec{I}}_{i-1\text{ terms}} 
    \otimes ~ \vec{s}^{\textup{x}/\textup{y}/\textup{z}} \otimes 
    \underbrace{\vec{I} \otimes \cdots \otimes \vec{I}}_{N-i\text{ terms}}
\end{align*}
where 
\begin{align*}
    \vec{s}^{\textup{x}} = \begin{bmatrix} 0 & 1 \\ 1 & 0 \end{bmatrix}
    && 
    \vec{s}^{\textup{y}} = \begin{bmatrix} 0 & -i \\ i & 0 \end{bmatrix}
    && 
    \vec{s}^{\textup{z}} = \begin{bmatrix} 1 & 0 \\ 0 & -1 \end{bmatrix}.
\end{align*}

\begin{figure}
    \centering
    \hfill\includegraphics[scale=.6]{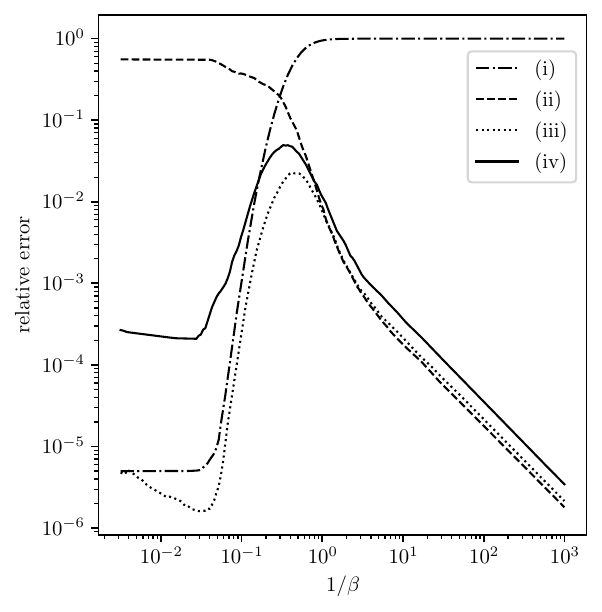}
    \hfill
    \includegraphics[scale=.6]{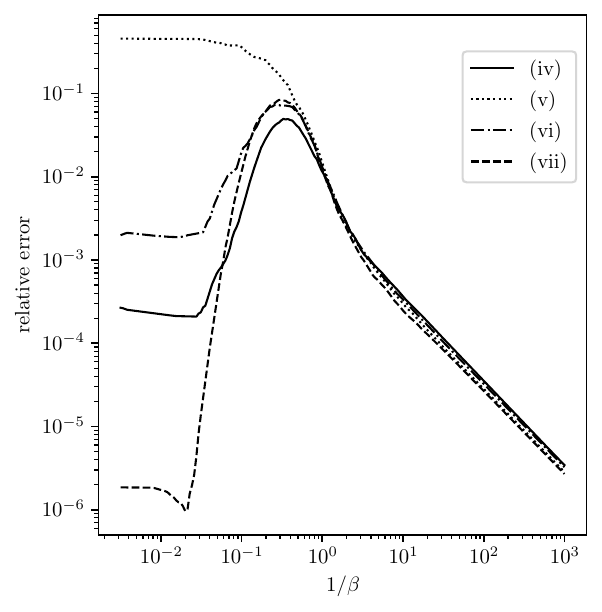}\hfill
    \caption{90th percentile of relative errors for \cref{alg:main} (left) and \cref{alg:restart} (right) used to estimate $\tr(\exp(-\beta\vec{A}))$ for a spin chain.
    Curve (i) corresponds to low-rank approximation only, curve (ii) corresponds to quadratic trace estimation only, curve (iii) corresponds to a combination whose total cost is the cost of (i)+(ii), and curve (iv) corresponds to a combination whose total cost is the same as each of (i) and (ii).
    Curve (v) corresponds to zero restarts, curve (vi) to two restarts, and curve (vii) to four restarts.
    Exact parameter choices are given in \cref{tab:spin_example}.
    \textbf{Note: the figure in the published version was generated without normalization, due to a clerical error.}
    }
    \label{fig:spin_example}
\end{figure}

We set $h=0.3$ and $N = 20$ so that $d = 2^{20} = 1048576$ and apply \cref{alg:main} or \cref{alg:restart} with several different choices of parameters to compute $Z(\beta)$ for a range of $\beta$.
For each parameter choice shown in \cref{tab:spin_example}, we run the algorithm independently 100 times and compare to the true value of $Z(\beta)$, which can be computed analytically for $s=1/2$ using the standard ``Bethe ansatz'' \cite{karabach1997introduction}.
In all cases, $n$ is fixed to be large enough that the matrix exponential is applied accurately. 
The 90-th percentile of the relative errors are reported in \cref{fig:spin_example}; i.e. errors were better than reported in the figure in 90\% of the trials.

\begin{table}
    \caption{Choices of parameters for \cref{fig:spin_example}. In the second-to-last column, an entry such as 1200+300=1500 indicates that 1200 matvecs were devoted to deflation and 300 to estimating the trace of the remainder, for a total of 1500.}
    \label{tab:spin_example}
    \centering
    \def\arraystretch{1.2}
    \begin{tabular}{ccccccrwr{.35cm}c} \toprule
        & $r$ & $b$ & $q$ & $m$ & $n$ & \multicolumn{2}{c}{\# matvecs} & $\operatorname{cols}(\vec{\bar{Q}}_{q+1})$ \\ \midrule
        (i) & 0 & 8 & 30 & 0 & 50 & $640+0=$&640 & 248 \\ 
        (ii) & 0 & 0 & 0 & 13 & 50 & $0+650=$&650 & 0 \\ 
        (iii) & 0 & 8 & 30 & 13 & 50 & $640+650=$&1290 & 248 \\ 
        (iv) & 0 & 4 & 30 & 6 & 50 & $320+300=$&650 & 124 \\ 
        (v) & 0 & 4 & 10 & 6 & 50 & $240+300=$&540 & 44 \\ 
        (vi) & 2 & 4 & 10 & 6 & 50 & $720+300=$&1020 & 44 \\ 
        (vii) & 4 & 4 & 10 & 6 & 50 & $1200+300=$&1500 & 44 \\ 
        \bottomrule
    \end{tabular}
\end{table}

The results for the first four parameter choices are shown in the left panel of \cref{fig:spin_example} and illustrate the benefit of our algorithm over approaches based solely on low-rank approximation or quadratic trace estimation.
Parameter choice (i) is a low-rank approximation type approach with $m=0$.
This produces a highly accurate approximation for large $\beta$, but a very inaccurate approximation for smaller $\beta$.
Parameter choice (ii) is a pure quadratic trace estimation approach with $b=q=0$.
This performs most accurately at small $\beta$ but less accurately at large $\beta$.
The next two parameter choices combine both approaches. 
Parameter choice (iii) uses the same non-zero values of $q$, $b$, and $m$ from choices (i) and (ii) and, unsurprisingly, performs better than the first two approaches.
Of course, the total number of matvecs is the sum of the first two approaches.
Parameter choice (iv) is around the same cost as the first two approaches, and while it performs somewhat worse than either of the original approaches at extreme values of $\beta$, it performs better than these approaches for intermediate values of $\beta$.
More importantly, the approach produces an approximation which is more uniformly good over the entire range of $\beta$ encountered.

We remark that the fact that the quality of the approximations improve as $\beta \to 0$ is due to the normalization.
If we do not use this normalization factor, then the error stagnates for $\beta$ small.
See \cite[\S 4.2]{epperly2023xtrace} for more experiments regarding the use of normalization.

The results for the final three parameter choices are shown in the right panel of \cref{fig:spin_example} and illustrate the potential effectiveness of restarting.
We take the filter polynomials to be the Chebyshev interpolants to $x\mapsto \exp(-\beta_0 x)$ on an interval containing the eigenvalues of $\vec{\bar{T}}_q$, where $\beta_0$ is the largest of the $\beta$ used in the experiment.
As expected, restarting allows a better low-rank approximation while limiting the storage and organization costs.
With sufficient restarts, we are able to attain an approximation of quality comparable or even better to the results of parameter choice (iv).

\subsection{Parameter selection for the inverse}

In this example, we explore how the block size $b$, Krylov depth $q$, and number of restart cycles $r$ impacts the quality of the basis $\vec{\bar{Q}}_{q+1}$.
Specifically, we compute 
\begin{equation}
\label{eqn:subspace_proj}
    \|(\vec{I} - \vec{\bar{Q}}_{q+1} \vec{\bar{Q}}_{q+1}^\T)f(\vec{A}) (\vec{I} - \vec{\bar{Q}}_{q+1} \vec{\bar{Q}}_{q+1}^\T)\|_{\F}/\|f(\vec{A})\|_{\F}
\end{equation} 
At least assuming $n$ is large enough that $f(\vec{A})$ is applied accurately, this quantity is directly proportional to the variance of the quadratic trace estimator used in the second stage of our algorithm.

In order to test our algorithm, we use $f(x) = 1/x$ and choose two spectra defined, for $i=1,2, \ldots, d$, by
\begin{equation}
\label{eqn:inv_spec}
f(\lambda_i^{\textup{slow}}) = 1 + \left(\frac{i-1}{d-1}\right)^2 (\kappa-1)
,\qquad
f(\lambda_i^{\textup{fast}}) = 1+ \left(\frac{i-1}{d-1}\right) (\kappa-1)   \rho^{d-i}.
\end{equation}
These respectively correspond to algebraic and geometric decay in the eigenvalues of the \emph{matrix function} $f(\vec{A})$.
We then use the first stage of \cref{alg:main,alg:restart} to compute $\vec{\bar{Q}}_{q+1}$ for a range of $q$, $b$, and $r$.

\begin{figure}
    \centering
    \hfill\includegraphics[scale=.7]{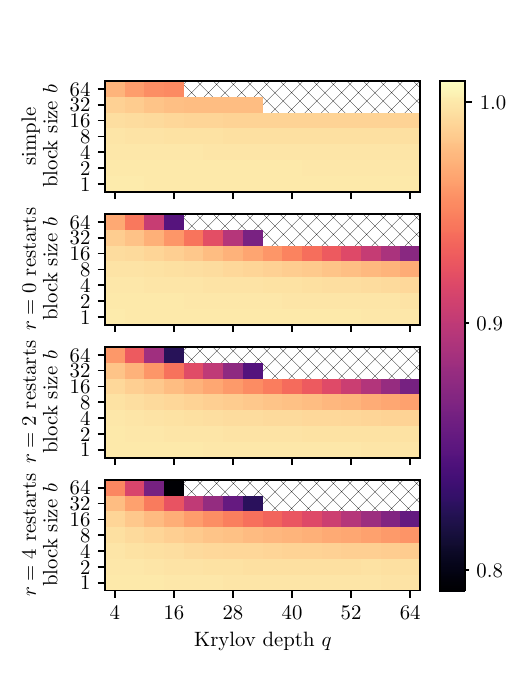}
    \hfill
    \includegraphics[scale=.7]{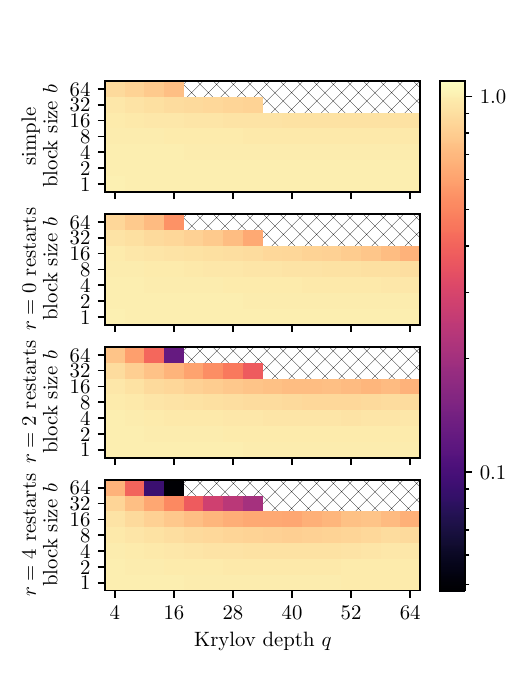}\hfill
    \caption{Quality of projection subspace $\vec{\bar{Q}}_{q+1}$ as in \cref{eqn:subspace_proj} (smaller is better) for $f(\vec{A})$ with algebraic decay \cref{eqn:inv_spec} with $n=2000$ and $\kappa = 1000$.
    Values for which $qb>1024$ were not computed. For reference, we also show the space obtained by \cref{alg:generalized_simple_func} as ``simple''.}
    \label{fig:proj_quality_slow}
\end{figure}

The results of our experiments are reported in \cref{fig:proj_quality_slow}.
As expected, when $q$, $b$, and $r$ are larger, the quality of the approximation improves.
In the case of slow (algebraic) decay, restarting is only mildly effective due to the fact that the are many eigenvalues with similar magnitude to the top eigenvalues.
Moreover, as expected, the reduction in the value of \cref{eqn:subspace_proj} is not substantial. 
On the other hand, in the case of fast (geometric) decay, restarting allows higher quality approximations.

In addition, we show the quality of projection used by \cref{alg:generalized_simple_func}, which is essentially \hpp{} with matrix-vector products with $f(\vec{A})$ computed via a black-box Krylov subspace method.
As expected, after $q$ becomes sufficiently large such that the Lanczos approximation to $f(\vec{A})\vec{\Omega}$ is reasonably accurate, this approach no longer improves with $q$. 
Moreover, for any fixed values of $q$ and $b$, the approximation is worse than than our algorithm, due to the fact that the projection space is a strict subspace of $\vec{\bar{Q}}_{q+1}$.

\subsection{Adaptive algorithm}
We test \cref{alg:adaptive} on several problems, one with synthetic data and two using matrices from the SuiteSparse Matrix Collection \cite{davis2011university}. On each of these problems we compare its performance to that of A-\hpp{} \cite{persson2022improved}. These algorithms differ only in their method for producing a set of vectors $\vec{\bar{Q}}_{q+1}$ for deflation and estimating $\tr(\vec{\bar{Q}}_{q+1}^\T f(\vec{A})\vec{\bar{Q}}_{q+1})$; in our implementation, the procedures for estimating the remainder (\hyperref[line:ada_remStart]{lines \ref{line:ada_remStart}}-\ref{line:ada_sampling} of \cref{alg:adaptive}) are identical. Thus for a fixed input $(\epsilon,\delta)$ we expect the algorithms to have similar accuracy. The primary differences will be the number of vectors devoted to deflation and the cost of the deflation step. 

\subsubsection{Estrada index}\label{sec:roget}
For our first test case we estimated the Estrada index $\tr(\exp(\vec{A}))$ of the Roget's Thesaurus graph, a $1022\times 1022$ matrix with 7297 nonzero elements after symmetrization. We ran the adaptive algorithms to a relative error of $2^{-p}$ for $2\leq p \leq 7$ with failure probability $\delta = 0.05$ and with 10 trials for each value of $p$. Matrix-vector products with $\exp(\vec{A})$ were estimated using $n=30$ steps of the Lanczos process, and the block-Lanczos routine in \cref{alg:adaptive} used block size $b=2$. 

\begin{figure}
    \centering
    \includegraphics[scale=.6]{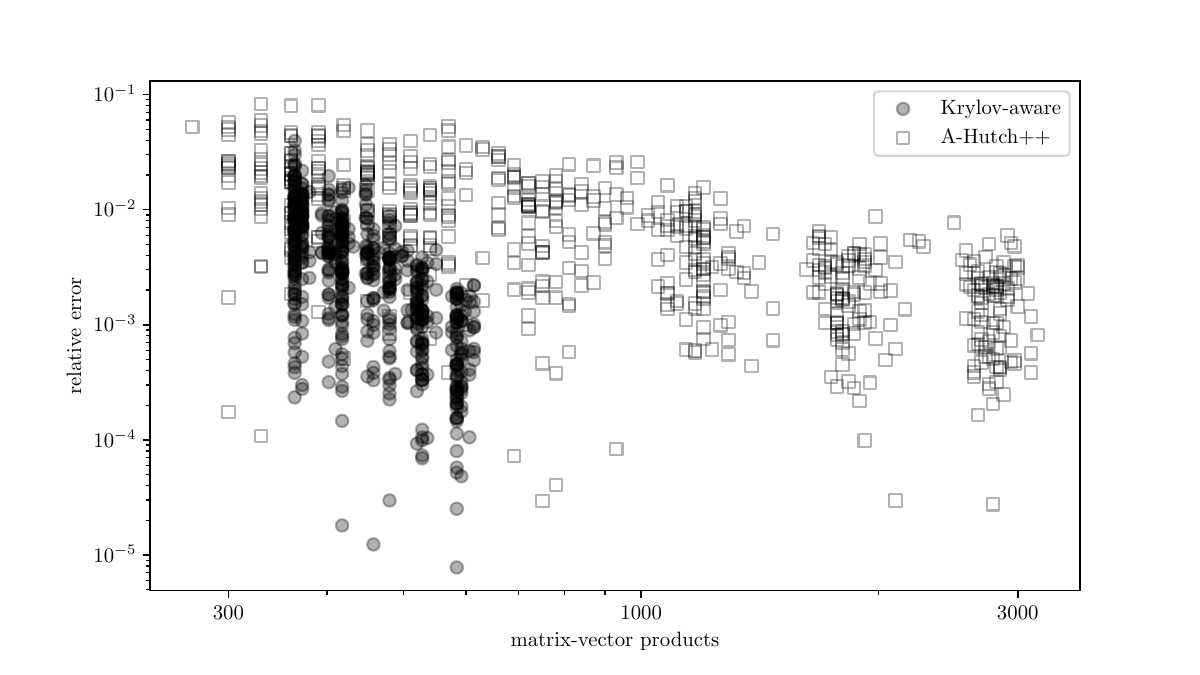}
    \caption{
    Costs associated with estimating $\tr(\exp(\vec{A}))$ for the Roget's Thesaurus graph.}
    \label{fig:roget_estrada}
\end{figure}

Results are shown in \cref{table:cost_estrada} and \cref{fig:roget_estrada}. From the table in particular, we make the following observations: 
\begin{itemize}[itemsep=0pt]
    \item \cref{alg:adaptive} used less than half as many matvecs as A-\hpp{} for $p=2$. As the desired relative error decreased, our algorithm's comparative advantage increased to nearly a factor of 7 for $p=7$. 
    \item \cref{alg:adaptive} used about 10 times as many vectors for deflation regardless of $\epsilon$. Even so, it required many fewer matvecs for the deflation step. 
    \item Because \cref{alg:adaptive} used more vectors for deflation, it did not need nearly as many samples to estimate the remainder. 
\end{itemize}
With the caveat that counting matvecs alone does not account for the cost of orthogonalizing and storing $\vec{\bar{Q}}_{q+1}$, it is apparent that drawing the deflation vectors from a block Krylov space and taking advantage of that structure has the potential to greatly reduce the cost of trace estimation problems. 

\begin{figure}
    \centering
    \includegraphics[scale=.6]{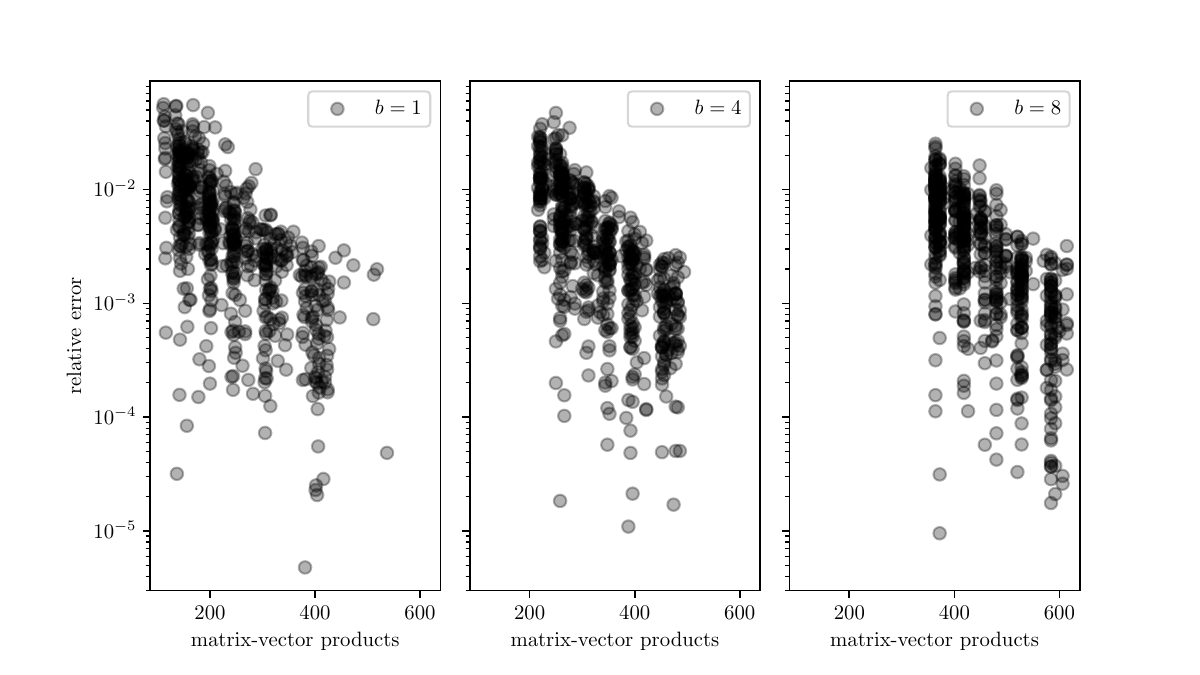}
    \caption{Costs associated with estimating $\tr(\exp(\vec{A}))$ for the Roget's Thesaurus graph with different block sizes.}
    \label{fig:roget_estrada_blocksize}
\end{figure}

We also ran the same set of experiments using block sizes $b\in \{1,4,8\}$. Results are shown in \cref{table:cost_estrada_blocksize} and \cref{fig:roget_estrada_blocksize}. Our algorithm (unsurprisingly) used fewer matvecs when the block size was smaller, although the difference became less pronounced for smaller error tolerances $\epsilon$. Smaller block sizes also tended to use fewer vectors for deflation. So at least to the extent that the number of matvecs is a reasonable proxy for the computational cost, we recommend using smaller block sizes over larger ones. 


\begin{table}
    \caption{Costs associated with estimating $\tr(\exp(\vec{A}))$ for the Roget's Thesaurus graph. Values are the average over 100 trials rounded to the nearest integer.}
    \label{table:cost_estrada}
    \centering
    \def\arraystretch{1.2}
    \begin{tabular}{cccrccrwr{.35cm}} \toprule
    & \multicolumn{3}{c}{\cref{alg:adaptive}} & \multicolumn{4}{c}{A-\hpp{}}\\ \cmidrule(lr){2-4} \cmidrule(lr){5-8}
    $p$ & $\operatorname{cols}(\vec{\bar{Q}}_{q+1})$ & $m$ & \multicolumn{1}{c}{\# matvecs} & $\operatorname{cols}(\vec{\bar{Q}}_{q+1})$ & $m$ & \multicolumn{2}{c}{\# matvecs}\\\midrule
2 & 72 & 2 & $304+60=364$ & 3 & 5 & $208 + 160 =$& 368 \\
3 & 80 & 2 & $312+74=386$ & 5 & 8 & $278 + 253 =$&  531 \\
4 & 96 & 3 & $328+93=421$ & 6 & 14 & $390 + 408 =$& 798 \\
5 & 130 & 4 & $362+106=469$ & 10 & 20 & $609 + 586 =$& 1195 \\
6 & 174 & 4 & $406+118=524$ & 15 & 31 & $916 + 936 =$&1851 \\
7 & 233 & 4 & $465+125=590$ & 24 & 45 & $1445 + 1354 =$& 2799 \\
    \bottomrule
\end{tabular}
\end{table}

\begin{table}
    \caption{Costs associated with estimating $\tr(\exp(\vec{A}))$ for the Roget's Thesaurus graph with \cref{alg:adaptive} using different block sizes. The number of vectors used for deflation is $\operatorname{cols}(\vec{\bar{Q}}_{q+1}) = (q+1)b$.  Values are the average over 100 trials rounded to the nearest integer.}
    \label{table:cost_estrada_blocksize}
    \centering
    \def\arraystretch{1.2}
    \begin{tabular}{cccccccccc} \toprule
    & \multicolumn{3}{c}{$b = 1$} & \multicolumn{3}{c}{$b = 4$}& \multicolumn{3}{c}{$b = 8$}\\ \cmidrule(lr){2-4} \cmidrule(lr){5-7} \cmidrule(lr){8-10}
    $p$ & $(q+1)b$ & $m$ & \# mv & $(q+1)b$ & $m$ & \# mv & $(q+1)b$ & $m$ & \# mv\\\midrule
    2  & 23  & 3 & 140  & 44   & 2 & 230  & 72  & 2 & 364 \\
    3  & 36  & 3 & 163  & 56   & 3 & 261  & 81  & 2 & 386 \\
    4  & 51  & 4 & 202  & 73   & 4 & 298  & 96  & 3 & 421 \\
    5  & 92  & 4 & 253 & 113   & 4 & 342 & 130  & 4 & 469 \\
    6 & 124  & 5 & 316 & 157   & 4 & 398 & 174  & 4 & 523 \\
    7 & 160  & 7 & 408 & 216   & 4 & 465 & 233  & 4 & 589 \\
    \bottomrule
\end{tabular}
\end{table}

\subsubsection{Log Determinant}
For our second test case we estimated the log determinant $\tr(\log(\vec{A}))$ of the matrix \texttt{thermomech\_TC}, a $102158\times 102158$ matrix with $711558$ nonzero elements. We ran the adaptive algorithms to a relative error of $2^{-p}$ for $6\leq p \leq 9$ with failure probability $\delta = 0.05$ and with 10 trials for each value of $p$. Matrix-vector products with $\log(\vec{A})$ were estimated using $n=35$ steps of the Lanczos process, and the block-Lanczos routine in \cref{alg:adaptive} used block size $b=2$. 

The \texttt{thermomech\_TC} matrix is very well-conditioned with $\kappa(\vec{A})\approx 70$, so not many samples are needed to get an accurate estimate and the variance reduction step is not particularly useful. Consequently, \cref{alg:adaptive} and A-\hpp{} have very similar performance. Both terminate the variance reduction step almost immediately; the only difference is that the overhead cost is lower for our method. 

\begin{table}
    \caption{Costs associated with estimating $\tr(\log(\vec{A}))$  for thermomechTC.  Values are the average over 50 trials rounded to the nearest integer.}
    \label{fig:cost_logdet}
    \centering
    \def\arraystretch{1.2}
    \begin{tabular}{cccrwr{.35cm}ccrwr{.35cm}} \toprule
    & \multicolumn{4}{c}{\cref{alg:adaptive}} & \multicolumn{4}{c}{A-\hpp{}}\\ \cmidrule(lr){2-5} \cmidrule(lr){6-9}
    $p$ & $\operatorname{cols}(\vec{\bar{Q}}_{q+1})$ & $m$ & \multicolumn{2}{c}{\# matvecs} & $\operatorname{cols}(\vec{\bar{Q}}_{q+1})$ & $m$ & \multicolumn{2}{c}{\# matvecs}\\\midrule
    5 & 6 &  3 & $74 +  105=$ &  179 & 3 &  3 & $210 +  105=$ &  315 \\
    6 & 6 &  4 & $74 +  140=$ &  214 & 3 &  4 & $210 +  140=$ &  350 \\
    7 & 6 &  8 & $74 +  280=$ &  354 & 3 &  8 & $210 +  280=$ &  490 \\
    8 & 6 & 19 & $74 +  665=$ &  739 & 3 & 19 & $210 +  665=$ &  875 \\
    9 & 6 & 55 & $74 + 1925=$ & 1999 & 3 & 55 & $210 + 1925=$ & 2135 \\
    \bottomrule
\end{tabular}
\end{table}

\begin{figure}
    \caption{Costs asscoiated with estimating $\tr(\log(\vec{A}))$ for thermomechTC.}
    \label{fig:logdet-thermo}
    \centering
    \includegraphics[scale=.6]{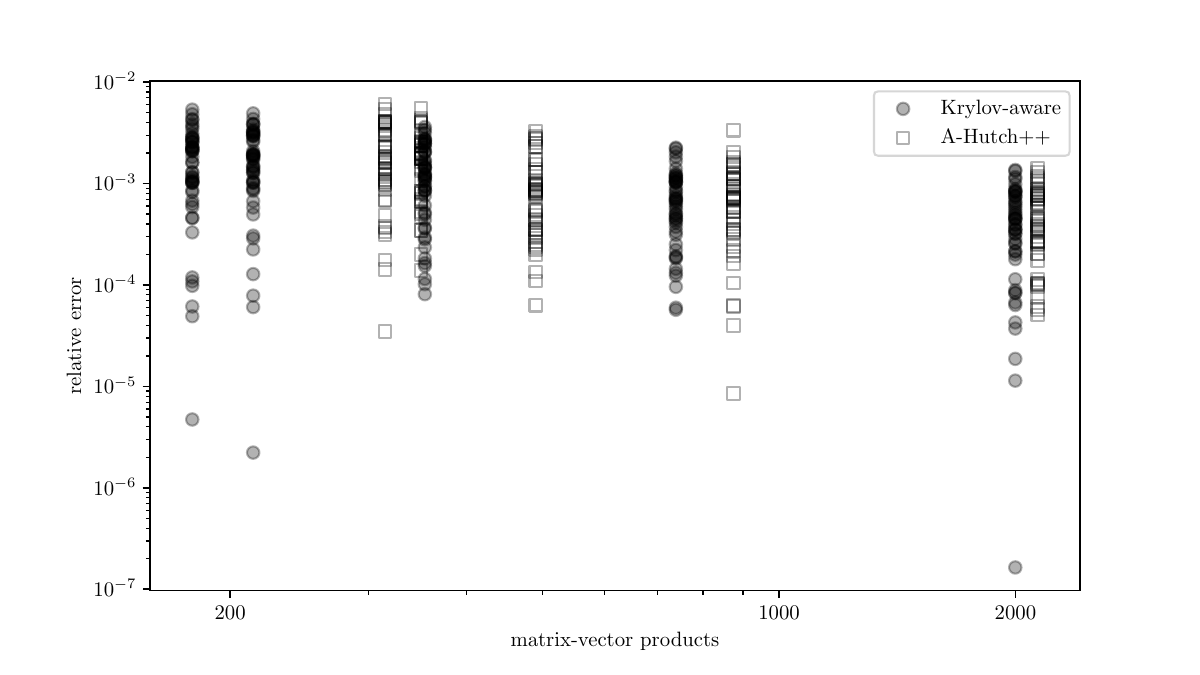}
\end{figure}

\subsubsection{Nuclear norm}
For our final test case we estimated $\tr(\vec{A}^{1/2})$ of a synthetic $2500\times 2500$ matrix $\vec{A} = \operatorname{diag}(1^{-c},2^{-c},\ldots,2500^{-c})$, with $c=1.5$. We ran the adaptive algorithms to a relative error of $2^{-p}$ for $2\leq p \leq 7$ with failure probability $\delta = 0.05$ and with 10 trials for each value of $p$. Matrix-vector products with $\vec{A}^{1/2}$ were estimated using $n=50$ steps of the Lanczos process, and the block-Lanczos routine in \cref{alg:adaptive} used block size $b=2$. 

\begin{figure}
    \centering
    \includegraphics[scale=.6]{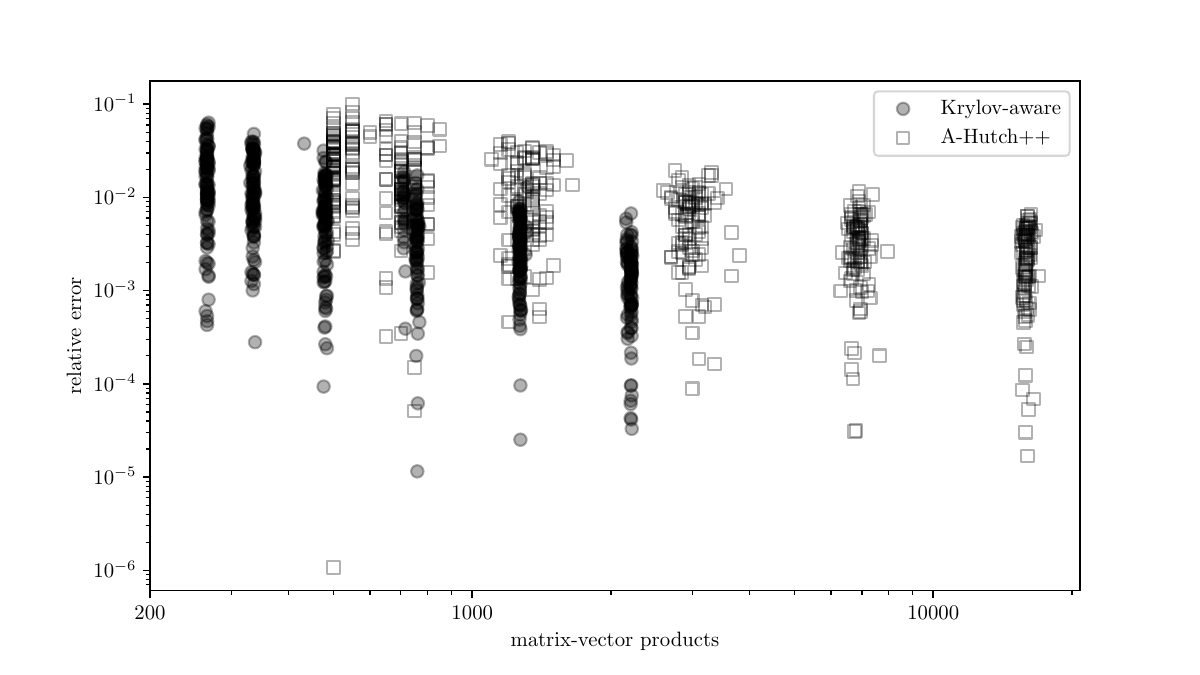}
    \caption{Estimating $\tr(\vec{A}^{1/2})$ for a synthetic matrix.}
    \label{fig:nuclear_norm}
\end{figure}

Results are shown in \cref{table:nuclear_norm} and \cref{fig:nuclear_norm}. The general pattern is similar to that of our experiments for the Roget's Thesaurus graph in \cref{sec:roget}: \cref{alg:adaptive} devoted more vectors to deflation (about 20 times as many as A-Hutch{\footnotesize ++}) and fewer toward estimating the remainder. It was more efficient overall, and the improvement increased from about a factor of 2 for $p=2$ to a factor of 7 for $p=7$. 

\begin{table}
    \caption{Costs associated with estimating $\tr(\vec{A}^{1/2})$ for a synthetic matrix. Values are the average over 100 trials rounded to the nearest integer.}
    \label{table:nuclear_norm}
    \centering
    \begin{tabular}{wc{.1cm}cwc{.25cm}rwr{.35cm}cwc{.25cm}rwr{.48cm}}  \toprule
    & \multicolumn{4}{c}{\cref{alg:adaptive}} & \multicolumn{4}{c}{A-\hpp{}}\\ \cmidrule(lr){2-5} \cmidrule(lr){6-9}
    $p$ & $\operatorname{cols}(\vec{\bar{Q}}_{q+1})$ & $m$ & \multicolumn{2}{c}{\# matvecs} & $\operatorname{cols}(\vec{\bar{Q}}_{q+1})$ & $m$ & \multicolumn{2}{c}{\# matvecs}\\\midrule
    2    &18     &3   &$116   +150=$&   266 &    3 &    4 &  $300   +216=$   & 516 \\
    3    &37     &4   &$135   +200=$&   335 &    3 &    8 &  $300   +419=$   & 719 \\
    4    &82     &6   &$180   +300=$&   479 &    4 &   19 &  $360   +962=$  & 1322 \\
    5   &162    &10   &$260   +488=$&   747 &    7 &   46 &  $726  +2286=$  & 3012 \\
    6   &320    &17   &$418   +852=$&  1270 &   18 &  102 & $1801  +5080=$  & 6881 \\
    7   &816    &26   &$914  +1286=$&  2199 &   44 &  232 & $4354 +11587=$ & 15941 \\
    \bottomrule
\end{tabular}
\end{table}

\section{Proofs} \label{section:proofs}
In this section, we provide proofs of the theoretical results stated above.

\begin{proof}[Proof of \cref{thm:poly_exact}]
Note that $\vec{\bar{Q}}_q\vec{\bar{Q}}_q^\T$ is the orthogonal projector onto $\mathcal{K}_q(\vec{A},\vec{Z})$, and for any $j\leq q-1$, $\vec{A}^j \vec{Z} \in \mathcal{K}_q(\vec{A},\vec{Z})$.
Thus, for any $k\leq q-1$, 
\begin{align*}
\vec{A}^k\vec{Z} 
= \vec{\bar{Q}}_q\vec{\bar{Q}}_q^\T \vec{A}^k\vec{Z}
&= \vec{\bar{Q}}_q\vec{\bar{Q}}_q^\T\vec{A} \vec{A}^{k-1} \vec{Z}
\\&= \vec{\bar{Q}}_q\vec{\bar{Q}}_q^\T\vec{A} \vec{\bar{Q}}_q\vec{\bar{Q}}_q^\T \vec{A}^{k-1} \vec{Z}
\\&= \vec{\bar{Q}}_q\vec{\bar{Q}}_q^\T\vec{A} \vec{\bar{Q}}_q\vec{\bar{Q}}_q^\T \vec{A} \vec{\bar{Q}}_q\vec{\bar{Q}}_q^\T \cdots \vec{\bar{Q}}_q\vec{\bar{Q}}_q^\T\vec{Z}.
\end{align*}
From \cref{eqn:krylov_recurrence} and the orthogonality of $\bar{\vec{Q}}_{q+1}$, we have that $\vec{\bar{Q}}_q^\T\vec{A} \vec{\bar{Q}}_q = \vec{T}_q$.
Thus since $\vec{\bar{Q}}_q^\T\vec{Z} = \vec{E}_1\vec{R}_1$, we find 
\[
\vec{A}^j \vec{Z} = \vec{\bar{Q}}_q [\vec{T}_q^j]_{:,1:b} \vec{R}_1.
\]
By linearity, \cref{eqn:mv_approx} is exact for any $f$ a polynomial of degree up to $q-1$.

As for \cref{eqn:quad_approx}, let $k\leq 2q-1$ be partitioned as $k = 2j+\ell$, where $j\leq q-1$ and $\ell \in \{0,1\}$. Recalling that $\vec{\bar{Q}}_q^\T \vec{A} \vec{\bar{Q}}_q = \vec{T}_q$, the previous result implies that
\begin{align*}
     \vec{Z}^\T \vec{A}^{k} \vec{Z} 
     &= (\vec{A}^{j}\vec{Z})^\T \vec{A}^\ell (\vec{A}^{j} \vec{Z})
     \\&= (\vec{\bar{Q}}_{q} (\vec{T}_{q})^{j} \vec{E}_1 \vec{R}_1)^\T 
     \vec{A}^\ell 
     (\vec{\bar{Q}}_{q} (\vec{T}_{q})^{j} \vec{E}_1 \vec{R}_1)
     \\&= ((\vec{T}_{q})^{j} \vec{E}_1 \vec{R}_1)^\T (\vec{\bar{Q}}_q^\T \vec{A}^\ell \vec{\bar{Q}}_q) ((\vec{T}_{q})^{j} \vec{E}_1 \vec{R}_1)
     \\&= \vec{R}_1^\T \vec{E}_1^\T (\vec{T}_{q})^{2j+\ell} \vec{E}_1 \vec{R}_1.
     \\&= \vec{R}_1^\T \vec{E}_1^\T (\vec{T}_{q})^{k} \vec{E}_1 \vec{R}_1.
\end{align*}
Again, by linearity, \cref{eqn:quad_approx} is exact for any $f$ a polynomial of degree up to $2q-1$.
\end{proof}

\begin{proof}[Proof of \cref{thm:func_rate}]
Given $\vec{A}$ and $\vec{Z}$, define for convenience
\[\mathsf{lan}_q(f) \equiv \vec{\bar{Q}}_q\left[f(\vec{T}_q)\right]_{:,1:b}\vec{R}_1\]
as the approximation \cref{eqn:mv_approx} to $f(\vec{A})\vec{Z}$. Let $p$ be any polynomial with $\deg(p)\leq q-1$. Then from the first part of \cref{thm:poly_exact}, the equality $\| \vec{R}_1 \| = \| \vec{Z} \|$, and the fact that the eigenvalues of $\vec{T}_q$ interlace those of $\vec{A}$, it holds that
\begin{align*}
    \| f(\vec{A})\vec{Z} - \mathsf{lan}_q(f) \|_2
    &\leq \| f(\vec{A})\vec{Z} - p(\vec{A})\vec{Z} \|_2 + \| \mathsf{lan}_q(p)-  \mathsf{lan}_q(f) \|_2
    \\& \leq \| \vec{Z} \|_2 \left( \| f(\vec{A}) - p(\vec{A}) \|_2 + \| f(\vec{T}_q) - p(\vec{T}_q) \|_2 \right) 
    \\&= \| \vec{Z} \|_2 \left( \max_{x\in\Lambda(\vec{A})} |f(x) - p(x)| + \max_{x\in\Lambda(\vec{T}_q)} |f(x) - p(x)| \right)
    \\& \leq 2 \| \vec{Z} \|_2 \max_{x\in [\lambda_{\textup{min}}, \lambda_{\textup{max}}]} | f(x) - p(x)|.
\end{align*}
The first result follows by optimizing over polynomials $p$. The second result is proved in a similar fashion using the second part of \cref{thm:poly_exact}.
\end{proof}

\begin{proof}[Proof of \cref{thm:hpp_func}]
By the triangle inequality and basic properties of the Frobenius norm,
\begin{align*}
    \| \vec{P}_{\vec{Q}} f(\vec{A}) \vec{P}_{\vec{Q}} 
    -  \vec{P}_{\vec{V}} f(\vec{A}) \vec{P}_{\vec{V}} \|_{\F}
    & \leq \| \vec{P}_{\vec{Q}} f(\vec{A}) \vec{P}_{\vec{Q}}  - 
    \vec{P}_{\vec{Q}} f(\vec{A}) \vec{P}_{\vec{V}} \|_{\F}
    \\&\hspace{3em}+ \| \vec{P}_{\vec{Q}} f(\vec{A}) \vec{P}_{\vec{V}} - \vec{P}_{\vec{V}} f(\vec{A}) \vec{P}_{\vec{V}} \|_{\F}
    \\&\leq 2 \|f(\vec{A})\|_2 \| \vec{P}_{\vec{Q}} - \vec{P}_{\vec{V}} \|_\F.
\end{align*}
Using the assumption that $\vec{V}$ and $\vec{Q}$ are of equal rank, \cite[Theorems 2.3, 2.4]{stewart1977perturbation} assert that
\[
    \| \vec{P}_{\vec{Q}}  - \vec{P}_{\vec{V}} \|_{2} 
    \leq \| (f(\vec{A}) \vec{\Omega})^\dagger \|_{2} \| f(\vec{A}) \vec{\Omega} - \mathsf{lan}_q(f) \|_2.
\]
Since $f(\vec{A})$ is square,
$\sigma_{\textup{min}}(f(\vec{A}) \vec{\Omega}) \geq  \sigma_{\textup{min}}(\vec{\Omega}) \sigma_{\textup{min}}( f(\vec{A}) )$. 
Thus, 
\[
\| (f(\vec{A})\vec{\Omega})^\dagger \| = \sigma_{\textup{min}}(f(\vec{A}) \vec{\Omega})^{-1} \leq \sigma_{\textup{min}}(\vec{\Omega})^{-1} \sigma_{\textup{min}}( f(\vec{A}) )^{-1}.
\]
Using this and the fact that $\operatorname{rank}(\vec{P}_{\vec{Q}}  - \vec{P}_{\vec{V}}) \leq 2b$,
\begin{align*}
    \| \vec{P}_{\vec{Q}} f(\vec{A}) \vec{P}_{\vec{Q}} 
    -  \vec{P}_{\vec{V}} f(\vec{A}) \vec{P}_{\vec{V}} \|_{\F}
    &\leq 
    2\| f(\vec{A}) \|_2 \|\vec{P}_{\vec{Q}}  - \vec{P}_{\vec{V}} \|_{\F}
    \\&\leq 2\sqrt{2b} \| f(\vec{A}) \|_2 \| \vec{P}_{\vec{Q}}  - \vec{P}_{\vec{V}} \|_{2}
    \\ &\leq 2\sqrt{2b} \| f(\vec{A}) \|_2 \| (f(\vec{A}) \vec{\Omega})^{\dagger} \|_2 \| f(\vec{A}) \vec{\Omega} - \mathsf{lan}_q(f) \|_2 
    \\&\leq 2\sqrt{2b}\, \frac{\sigma_{\textup{max}}(f(\vec{A}))}{\sigma_{\textup{min}}(f(\vec{A}))} \frac{\sigma_{\textup{max}}(\vec{\Omega})}{\sigma_{\textup{min}}(\vec{\Omega})} \Delta.
\end{align*}
\end{proof}

\begin{proof}[Proof of \cref{thm:krylov_quad_approx}]
The proof is essentially identical to the proof of \cref{thm:poly_exact} after relabeling $\vec{\bar{Q}}_{q} \to \vec{\bar{Q}}_{q+n+1}$, $\vec{Z} \to \vec{\bar{Q}}_{q}$, and $q\to n$. Since $\vec{\bar{Q}}_{q}$ has orthogonal columns, the analogue of $\vec{R}_1$ from \cref{thm:poly_exact} is just the identity.

\end{proof}

\begin{proof}[Proof of \cref{thm:main}]
Define
\[
    \mathsf{est} \equiv \tr(\vec{\hat{Q}}^\T f(\vec{A}) \vec{\hat{Q}})  + \frac{d-\hat{b}}{m} \sum_{i=1}^{m} \frac{\vec{y}_i^\T f(\vec{A})\vec{y}_i}{\vec{y}_i^\T \vec{y}_i} .
\]
Recall that $\vec{y}_i = (\vec{I} - \vec{\hat{Q}}\vec{\hat{Q}}^\T) \bm{\uppsi}_i$, where $\bm{\uppsi}_i$ is a Gaussian vector, so $\vec{y}_i/\| \vec{y}_i \|_2$ has a uniform distribution on the unit hypersphere defined on the complement of the column span of $\hat{\vec{Q}}$ which has dimension $d-\hat{b}$.
Thus $\EE[ \mathsf{est} ] = \tr(f(\vec{A}))$, and by the law of total variance and the variance formula for trace estimation using random vectors from the real hypersphere given in \cite{girard1987algorithme},
\begin{align*}
    \VV\big[\mathsf{est}\big] &= \EE\big[\VV[\mathsf{est} | \vec{\Omega}]\big] + \VV\big[\EE[\mathsf{est}|\vec{\Omega}]\big]\\
    &= \EE\left[\frac{2(d-\hat{b})}{m(d-\hat{b}+2)}\left(  \|\vec{F}\|_\F^2 - \frac{\tr(\vec{F})^2}{d-\hat{b}} \right)\right] + \underbrace{\VV\big[\tr(f(\vec{A}))\big]}_{=\,0}\\
    &= \frac{2(d-\hat{b})}{m(d-\hat{b}+2)}\left(  \EE\big[\|\vec{F}\|_\F^2\big] - \frac{\EE\big[\tr(\vec{F})^2\big]}{d-\hat{b}} \right).
\end{align*}
\Cref{alg:generalized_simple_func,alg:main} compute the quantities 
\[
    t_{\text{defl}} = \tr\left( \left[f(\vec{T}_{q+n})\right]_{1:\hat{b},1:\hat{b}} \right)
    ,\qquad
    t_{\text{rem}} =  \frac{d-\hat{b}}{m}\sum_{i=1}^{m}[f(\vec{T}_{n}^{(i)})]_{1,1},
\]
which respectively approximate the first and second terms of $\mathsf{est}$.
Since $\vec{\hat{Q}}$ has orthonormal columns, $\| \vec{\hat{Q}} \|_2^2 = 1$.
Thus, by the definition of $\Delta$,
\[
    \left\| \vec{\hat{Q}}^\T f(\vec{A}) \vec{\hat{Q}} - \left[f(\vec{T}_{q+n})\right]_{1:\hat{b},1:\hat{b}} \right\|_2 \leq \Delta,
\qquad
    \left| \frac{\vec{y}_i^\T f(\vec{A})\vec{y}_i}{\vec{y}_i^\T \vec{y}_i} - [f(\vec{T}_{n}^{(i)})]_{1,1} \right| \leq  \Delta.
\]
Using the fact that $|\!\tr(\vec{E})| \leq \hat{b}\, \|\vec{E}\|_2$ for any $\vec{E} \in \mathbb{R}^{\hat{b}\times \hat{b}}$,
\[
    \left| \mathsf{est} - (t_{\textup{defl}} +t_{\textup{rem}})\right|
    \leq \hat{b}\, \Delta + (d-\hat{b})\Delta
    = d\, \Delta.
\]
Taking expected values gives the bound for the expectation.

We can bound the variance by
\begin{equation*}
    \VV\big[ t_{\textup{defl}} +t_{\textup{rem}}   \big] \leq \left(\sqrt{\VV\big[ \mathsf{est}  \big]} + \sqrt{\VV\big[ \mathsf{est} - (t_{\textup{defl}} +t_{\textup{rem}}) \big]} \right)^2.
\end{equation*}
Since $\VV\big[ \mathsf{est} - (t_{\textup{defl}} +t_{\textup{rem}}) \big] \leq \EE \big[|\mathsf{est} - (t_{\textup{defl}} +t_{\textup{rem}})|^2 \big]\leq \EE\big[(d\, \Delta)^2\big] = d^2\,\EE\big[\Delta^2\big]$, we get the variance bound.

Finally, if the same $\vec{\Omega}$ is used in both algorithms, the column span of $\vec{Q}$ is that of $\vec{\bar{Q}}_{q}\left[f(\vec{T}_{q})\right]_{:,1:b}\vec{R}_1$, which is clearly contained in $\vec{\bar{Q}}_{q+1}$.
Thus, $\vec{P}_{\vec{\bar{Q}}_q} = \vec{P}_{\vec{\bar{Q}}_q} \vec{P}_{\vec{Q}}$, so 
basic properties of the Frobenius norm, and the fact that $\| \vec{P}_{\vec{\bar{Q}}_{q+1}}\|_2 \leq 1$, imply that
\[
    \| \vec{P}_{\vec{\bar{Q}}_{q+1}} f(\vec{A}) \vec{P}_{\vec{\bar{Q}}_{q+1}} \|_\F^2
    = \| \vec{P}_{\vec{\bar{Q}}_{q+1}}\vec{P}_{\vec{Q}} f(\vec{A}) \vec{P}_{\vec{Q}}\vec{P}_{\vec{\bar{Q}}_{q+1}} \|_\F^2
    \leq \| \vec{P}_{\vec{Q}} f(\vec{A}) \vec{P}_{\vec{Q}} \|_\F^2. 
\]
Taking expectations on both sides preserves this inequality.
\end{proof}

\section{Conclusions and future work}

From our analysis and experiments it is clear that exploiting the structure of block Krylov subspaces can significantly reduce the number of matrix-vector products required for randomized trace estimation. The effect on the computational time is less clear, and depends both on the cost of maintaining the orthonormal basis $\vec{\bar{Q}}_q$ and on how efficiently matvecs with $\vec{A}$ can be computed in parallel. A high-quality practical implementation for large-scale problems will likely require further study in order to more effectively balance the true costs of the algorithm. For such a setting, we believe that the restarted variant (\cref{alg:restart}), in particular, merits further study.

It is also worth considering the role of the approximation degree $n$ in further detail. 
For example, if the deflation step has proved useful and significantly reduced the Frobenius norm of the remainder, then it may be possible to estimate the trace of the remainder using a smaller value of $n$. We are not aware of any formal results on this topic. Informally, the paper \cite{hallman2022multilevel} proposes a trace estimator that uses multiple different degrees $n$, and in doing so hedges against overestimating $n$. 

\section{Acknowledgements}

The authors thank the referees and editor for their comments which improved the presentation of the paper. We also thank David Persson for helpful discussions.

\bibliography{references}
\bibliographystyle{abbrv}

\end{document}